\documentclass[12pt]{article}

\usepackage{amsmath}
\usepackage{amssymb}
\usepackage{color}
\usepackage{graphicx}
\usepackage{cite}
\usepackage{epstopdf}

\newtheorem{propo}{Proposition}
\newtheorem{rem}{Remark}

\newcommand{\rset}{{\mathbb{R}}}

\title{Convergence acceleration of Kaczmarz's method}
\author{Claude Brezinski\thanks{Laboratoire Paul Painlev\'e, UMR CNRS 8524, UFR
de Math\'ematiques Pures et Appliqu\'ees, Universit\'e des Sciences
et Technologies de Lille,
59655--Villeneuve d'Ascq cedex,
France, E--mail: {\tt Claude.Brezinski@univ-lille1.fr}.}
\and Michela Redivo--Zaglia\thanks{Universit\`a degli Studi di Padova,
Dipartimento di Matematica,
Via Trieste 63, 35121--Padova,
Italy. E--mail: {\tt Michela.RedivoZaglia@unipd.it}.
}
}

\begin{document}

\maketitle

\centerline{\it This paper is dedicated to the memory of Bernard Germain--Bonne (1940--2012)}

\begin{abstract}
The method of alternation projections (MAP) is an iterative procedure for finding the projection of a
point on the intersection of closed subspaces of an Hilbert space. The convergence of this method is
usually slow, and several methods for its acceleration have already been proposed. In this work, we consider
a special MAP, namely Kaczmarz' method for solving systems of linear equations.
The convergence of this method is discussed. After giving its matrix formulation
and its projection properties, we consider several procedures for accelerating its convergence.
They are based on sequence transformations whose kernels contain sequences of the same form as
the sequence of vectors generated by Kaczmarz' method. Acceleration can be achieved either directly,
that is without modifying the sequence obtained by the method, or by restarting it
from the vector obtained by acceleration. Numerical examples show the effectiveness
of both procedures.
\end{abstract}

\section{Introduction}

Let $Q_i$ denote the orthogonal projection on a closed subspace $M_i$ of an Hilbert space $H$, and let
$Q_M$ be the composition of the $r$ projecting operators $Q_i$, that is $Q_M=Q_r\cdots Q_1$. Let
$P_M$ be the projection on $M$, the intersection of the subspaces $M_i$.

We are looking for the projection of a given point ${\mathbf{\bar x}}$
on $M$. It holds
$$\lim_{n \to \infty}Q_M^n{\mathbf{\bar x}}=P_M{\mathbf{\bar x}}, \quad \forall {\mathbf{\bar x}} \in H.$$

The {\it method of alternating projections} (MAP) consists in the iterations
$${\mathbf{x}}_{n+1}=Q_M{\mathbf{x}}_n, \quad n=0,1,\ldots, \quad {\mathbf{x}}_0={\mathbf{\bar x}}.$$
This method often converges quite slowly and it needs to be accelerated \cite{marcos}. For this purpose,
the projection operator $Q_M$ can be replaced by another (non necessarily linear) operator $T$ which can
depend on $n$.

\vskip 2mm

There are two ways of using $T$
\begin{enumerate}
\item Keep the sequence $({\mathbf{x}}_n)$ as given by MAP
unchanged, and consider the new sequence $({\mathbf{y}}_n)$ built
by ${\mathbf{y}}_n=T{\mathbf{x}}_n$. \item Iterate the operator
$T$, that is consider the new iterations
${\mathbf{x}}_{n+1}=T{\mathbf{x}}_n$.
\end{enumerate}

\vskip 2mm

There exist many choices for $T$ (the norms are the Euclidean ones). For example, the choice
$$T{\mathbf{x}}_n={\mathbf{x}}_n-t_n(Q_M{\mathbf{x}}_n-{\mathbf{x}}_n),$$
with
$t_n=({\mathbf{x}}_n,{\mathbf{x}}_n-Q_M{\mathbf{x}}_n)/\|{\mathbf{x}}_n-Q_M{\mathbf{x}}_n\|^2$
was proposed in \cite[p. 44]{marcos}, while
 $t_n=({\mathbf{r}}_n,{\mathbf{v}}_n)/({\mathbf{v}}_n,{\mathbf{v}}_n)$
where ${\mathbf{r}}_n=Q_M{\mathbf{x}}_n-{\mathbf{x}}_n$ and
${\mathbf{v}}_n=Q_M^2{\mathbf{x}}_n-2Q_M{\mathbf{x}}_n+{\mathbf{x}}_n$
was discussed in \cite{gea}. More choices are presented in \cite{andrzej}.

\vskip 2mm

Kaczmarz's method \cite{kacz} for solving a system of linear equations was proposed on 1937.
Later, it was rediscovered by Gordon et al. \cite{art}, and applied in medical imagining. They called it
{\sc ART} ({\it Algebraic Reconstruction Technique}), and its original version, or some variants, continue to be used for tomographic imaging.
It is a particular case of row projection methods which received much attention (see,
for example, \cite{proj,ake,sameh,tana,andrzej}), and it enters into the framework of MAP.
It is also well suited for parallel
computations and large--scale problems, since each step only requires one row of the
matrix $A$ (or several rows simultaneously in its block version), and no matrix--vector products.
For an impressive list of publications on Kaczmarz's method, see \cite{list}.
Kaczmarz's method is also often used for solving an overdetermined consistent
$M \times N$ linear system with $M \geq N$, but, in this paper, we restrict ourselves to the case of a square regular
system. It is easy to extend the algorithms and the theory by properly replacing
$N$ by $M$.

A drawback of Kaczmarz's method, as in general of all projection iterative methods, is its often slow
convergence. Thanks to its matrix analysis, we will
be able to show how its convergence can be accelerated by some particular choices of the operator $T$ corresponding to the two ways of using it described above.

For definitions and properties of projections, see, for example, \cite{proj}.
In the sequel, the same notation will be used for a matrix and the projection it represents.

\section{Kaczmarz's method}

We consider a $N \times N$ linear system
$A{\mathbf{x}}={\mathbf{b}}$.
One single step of Kaczmarz's method consists in
\begin{equation}
{\mathbf{p}}_{n+1}=
{\mathbf{p}}_n+\frac{({\mathbf{b}}-A{\mathbf{p}}_n,{\mathbf{e}}_i)}{(A^T{\mathbf{e}}_i,A^T{\mathbf{e}}_i)}
A^T{\mathbf{e}}_i,
\label{kac}
\end{equation}
where ${\mathbf{e}}_i$ is the $i$th vector of the canonical basis of $\rset^N$.
There exist several strategies for choosing the index $i$ at each step. The most common one is
$i=n~ (\mbox{mod.}~ N)+1$. In this case, the method is called the {\it cyclic Kaczmarz's method} or, simply, the
{\it  Kaczmarz's method}, since it was originally proposed by Kaczmarz under this form \cite{kacz}.
It corresponds to a restarting from the result obtained after $N$ single steps, that is, in other words, to a renumbering of them or, again in other words, to the extraction of a subsequence.
Thus, one iteration of Kaczmarz's method consists in a complete cycle of steps in their natural order, that is

\begin{equation}
\label{cycle}
\left.
\begin{array}{rcl}
&&{\mathbf{p}}_0={\mathbf{x}}_n\\
&&{\mathbf{p}}_i={\mathbf{p}}_{i-1}+\displaystyle{\frac{({\mathbf{b}}-A{\mathbf{p}}_{i-1},
{\mathbf{e}}_i)}{(A^T{\mathbf{e}}_i,A^T{\mathbf{e}}_i)}}A^T{\mathbf{e}}_i, \qquad i=1,\ldots,N\\
&&{\mathbf{x}}_{n+1}={\mathbf{p}}_N
\end{array}
\right\}
\end{equation}

Denote by ${\mathbf a}_i=A^T{\mathbf{e}}_i$ the column vector formed by the $i$th row of $A$.
Thus, the computation of each vector ${\mathbf{p}}_i$ in
\eqref{cycle} does not require any matrix--vector product since
$({\mathbf{b}}-A{\mathbf{p}}_{i-1},{\mathbf{e}}_i)=({\mathbf{b}},
{\mathbf{e}}_i)-({\mathbf{p}}_{i-1},A^T{\mathbf{e}}_i)$. Thus, if
we denote by $b_i$ the $i$th component of the
right hand side ${\mathbf{b}}$, then the computation of
${\mathbf{p}}_{i}$ is simply given by

\begin{equation}
\label{mmm}
{\mathbf{p}}_i={\mathbf{p}}_{i-1}+\displaystyle{\frac{b_i -
({\mathbf{p}}_{i-1},{\mathbf{a}}_{i})}{\|{\mathbf{a}}_{i}\|^2}}\,{\mathbf{a}}_i.
\end{equation}

This remark is one of the main advantages of Kaczmarz' method, and
it allows an easy parallel implementation.

\vskip 2mm

Let
$M_i=\{{\mathbf{y}}~|~({\mathbf{b}}-A{\mathbf{y}},{\mathbf{e}}_i)=0\}$.
Then ${\mathbf{p}}_i$ is the oblique projection of
${\mathbf{p}}_{i-1}$ on $M_i$ along $A^T{\mathbf{e}}_i$.

\vskip 2mm

Moreover, setting
$$\lambda_i=\frac{({\mathbf{b}}-A{\mathbf{p}}_{i-1},{\mathbf{e}}_i)}{(A^T{\mathbf{e}}_i,A^T{\mathbf{e}}_i)},$$
one iteration of Kaczmarz's method  \eqref{cycle} writes
$${\mathbf{x}}_{n+1}={\mathbf{x}}_n+A^T\Lambda_n, \qquad \Lambda_n=(\lambda_1,\ldots,\lambda_N)^T \in \mathbb{R}^N.$$

\begin{rem}
~~\\
Let ${\mathbf{y}}_n$ be the iterates obtained by applying the
Gauss--Seidel method to the system $AA^T{\mathbf{y}}={\mathbf{b}}$. Then
${\mathbf{x}}=A^T{\mathbf{y}}$ and ${\mathbf{x}}_n=A^T{\mathbf{y}}_n$ \cite{altm}.
\end{rem}

Let us now analyze each step of \eqref{cycle}.
Inside one iteration, we have the following orthogonality properties for
$i=1,\ldots,N$ (see \cite{dura,gast})

\begin{eqnarray*}
&&({\mathbf{b}}-A{\mathbf{p}}_i,{\mathbf{e}}_i)=0, \\
&&({\mathbf{p}}_i-{\mathbf{p}}_{i-1},{\mathbf{x}}-{\mathbf{p}}_i)=0,
\end{eqnarray*}
and it follows
\begin{eqnarray*}
\|{\mathbf{x}}-{\mathbf{p}}_{i-1}\|^2&=&\|{\mathbf{x}}-{\mathbf{p}}_i\|^2+2\lambda_i({\mathbf{x}}-{\mathbf{p}}_i,A^T
{\mathbf{e}}_i)+\lambda_i^2(A^T{\mathbf{e}}_i,A^T{\mathbf{e}}_i), \\
&=&\|{\mathbf{x}}-{\mathbf{p}}_i\|^2+2\lambda_i({\mathbf{b}}-A{\mathbf{p}}_i,{\mathbf{e}}_i)+\|{\mathbf{p}}_i
-{\mathbf{p}}_{i-1}\|^2,
\end{eqnarray*}
since $({\mathbf{b}}-A{\mathbf{p}}_i,{\mathbf{e}}_i)=0$. Thus
$$\|{\mathbf{x}}-{\mathbf{p}}_{i-1}\|^2=\|{\mathbf{x}}-{\mathbf{p}}_i\|^2+\|{\mathbf{p}}_i-{\mathbf{p}}_{i-1}\|^2,$$
which shows that $\|{\mathbf{x}}-{\mathbf{p}}_i\| \leq
\|{\mathbf{x}}-{\mathbf{p}}_{i-1}\|$. Therefore, $\|{\mathbf{x}}-{\mathbf{x}}_{n+1}\| \leq
\|{\mathbf{x}}-{\mathbf{x}}_n\|$. Thus, as proved in \cite{gast}, these inequalities are strict, and Kaczmarz's method is always converging to the solution of the system.

Obviously $({\mathbf{b}}-A{\mathbf{p}}_{i-1},{\mathbf{e}}_i)^2
\leq
({\mathbf{b}}-A{\mathbf{p}}_{i-1},{\mathbf{b}}-A{\mathbf{p}}_{i-1})$.
In the case where
$({\mathbf{b}}-A{\mathbf{p}}_{i-1},{\mathbf{e}}_i)^2 \geq
({\mathbf{b}}-A{\mathbf{p}}_{i-1},{\mathbf{b}}-A{\mathbf{p}}_{i-1})/N$,
one can prove, by an analysis similar to what is done in \cite[p.
122]{dura}, that the following result holds

\begin{propo}
~~\\
If $({\mathbf{b}}-A{\mathbf{p}}_{i-1},{\mathbf{e}}_i)^2 \geq
({\mathbf{b}}-A{\mathbf{p}}_{i-1},{\mathbf{b}}-A{\mathbf{p}}_{i-1})/N$,
then
$$\|{\mathbf{x}}-{\mathbf{p}}_{i}\|^2 \leq \left(1-\frac{1}{N \kappa(AA^T)}\right)\|{\mathbf{x}}-{\mathbf{p}}_{i-1}\|^2.$$
\end{propo}

Then, in the general case, this coefficient could be an approximation of the convergence factor of the method.
This result is quite similar to a corresponding result proved in \cite{house} for Gastinel's method \cite{gast1} and other methods, or to an extension (Thm. 4.27 in \cite{gala}) of another result given in \cite{house}.

\subsection{Matrix interpretation}
\label{MatInt}
Following \cite{gast}, where it seems that it first
appeared, let us give the matrix interpretation of Kaczmarz's
method \eqref{cycle} (see also \cite{dura,sameh,tana}).

\vskip 2mm

We set
\begin{eqnarray*}
{\boldsymbol\alpha}_i&=&\displaystyle \frac{A^T{\mathbf{e}}_i}{\|A^T{\mathbf{e}}_i\|^2}=\frac{{\mathbf a}_i}
{\|{\mathbf a}_i\|^2},\\
P_i&=&I-A{\boldsymbol\alpha}_i {\mathbf{e}}_i^T, \\
Q_i&=&A^{-1}P_iA, \\
{\boldsymbol\varrho}_i&=&{\mathbf{b}}-A{\mathbf{p}}_i.
\end{eqnarray*}

We have
$$Q_i=A^{-1}P_iA=I-\frac{A^T{\mathbf{e}}_i{\mathbf{e}}_i^TA}{\|A^T{\mathbf{e}}_i\|^2}=I-\frac{{\mathbf a}_i
{\mathbf a}_i^T}{\|{\mathbf a}_i\|^2}=I-{\boldsymbol\alpha}_i {\mathbf a}_i^T.$$
The matrix $P_i$ represents the oblique projection on
${\mathbf{e}}_i^\perp$ along $AA^T{\mathbf{e}}_i$, while $Q_i$ is
the rank $N-1$ orthogonal projection on
$(A^T{\mathbf{e}}_i)^\perp$ along $A^T{\mathbf{e}}_i$. Thus, for
any vector ${\mathbf{y}}$, $(P_i{\mathbf{y}},{\mathbf{e}}_i)=0$
and $(Q_i{\mathbf{y}},A^T{\mathbf{e}}_i)=0$.

\vskip 2mm

Inside one iteration, we have, for $i=1,\ldots,N$,
\begin{equation}
\label{rhoi}
\left.
\begin{array}{rcl}
{\mathbf{p}}_i&=&Q_i{\mathbf{p}}_{i-1}+({\mathbf{b}},{\mathbf{e}}_i){\boldsymbol\alpha}_i=Q_i{\mathbf{p}}_{i-1}+(I-Q_i)
{\mathbf{x}}\\
{\mathbf{x}}-{\mathbf{p}}_i&=&{\mathbf{x}}-{\mathbf{p}}_{i-1}-({\boldsymbol\varrho}_{i-1},{\mathbf{e}}_i)
{\boldsymbol\alpha}_i=Q_i({\mathbf{x}}-{\mathbf{p}}_{i-1})\\
{\boldsymbol\varrho}_i&=&{\boldsymbol\varrho}_{i-1}-({\boldsymbol\varrho}_{i-1},{\mathbf{e}}_i)A{\boldsymbol\alpha}_i
=P_i{\boldsymbol\varrho}_{i-1}.
\end{array}
\right\}
\end{equation}
The first relation can be written
$${\mathbf{p}}_i=Q_i{\mathbf{p}}_{i-1}+A^{-1}(I-P_i){\mathbf{b}}.$$
Notice that
$A^{-1}(I-P_i){\mathbf{b}}=({\mathbf{b}},{\mathbf{e}}_i){\boldsymbol\alpha}_i$.

\vskip 2mm

The second relation in \eqref{rhoi} can also be written as
$${\mathbf{x}}-{\mathbf{p}}_i={\mathbf{x}}-{\mathbf{p}}_{i-1}-{\boldsymbol\alpha}_i {\mathbf{e}}_i^T {\boldsymbol\varrho}_{i-1},$$
and, replacing ${\boldsymbol\alpha}_i$ by its expression, it
follows
\begin{eqnarray*}
\|{\mathbf{x}}-{\mathbf{p}}_i\|^2&=&\|{\mathbf{x}}-{\mathbf{p}}_{i-1}\|^2-2\frac{({\boldsymbol\varrho}_{i-1},{\mathbf{e}}_i)}
{\|A^T{\mathbf{e}}_i\|^2}({\mathbf{x}}-{\mathbf{p}}_{i-1},A^T{\mathbf{e}}_i)+\frac{({\boldsymbol\varrho}_{i-1},
{\mathbf{e}}_i)^2}{\|A^T{\mathbf{e}}_i\|^2}\\
&=&\|{\mathbf{x}}-{\mathbf{p}}_{i-1}\|^2-\frac{({\boldsymbol\varrho}_{i-1},{\mathbf{e}}_i)^2}{\|A^T{\mathbf{e}}_i\|^2}.
\end{eqnarray*}

Summing up this identity for $i=1,\ldots,N$, we obtain an expression for the gain of one iteration of Kaczmarz's method

$$\|{\mathbf{x}}-{\mathbf{x}}_{n+1}\|^2=\|{\mathbf{x}}-{\mathbf{x}}_{n}\|^2- \sum_{i=1}^N \frac{({\boldsymbol\varrho}_{i-1},{\mathbf{e}}_i)^2}{\|A^T{\mathbf{e}}_i\|^2}.$$

Setting
\begin{eqnarray*}
P&=&P_N\cdots P_1\\
Q&=&Q_N\cdots Q_1=A^{-1}P_N\cdots P_1A=A^{-1}PA\\
{\mathbf{r}}_n&=&{\mathbf{b}}-A{\mathbf{x}}_n,
\end{eqnarray*}
it holds
\begin{equation}
\label{Qc}
\left.
\begin{array}{rcl}
{\mathbf{x}}-{\mathbf{x}}_{n+1}&=&Q({\mathbf{x}}-{\mathbf{x}}_n)\\
{\mathbf{r}}_{n+1}&=&P{\mathbf{r}}_n.
\end{array}
\right\}
\end{equation}
The matrix $Q$ represents an orthogonal projection, but not $P$. It obviously follows
\begin{equation}
\label{pqk}
{\mathbf{x}}-{\mathbf{x}}_n=Q^n({\mathbf{x}}-{\mathbf{x}}_0),
\quad {\mathbf{r}}_n=P^n{\mathbf{r}}_0,\quad n=0,1,\ldots,
\end{equation}
and we  have
\begin{equation}
\label{Qcc}
\begin{array}{rcl}
{\mathbf{x}}_{n+1}&=&Q{\mathbf{x}}_n+A^{-1}(I-P){\mathbf{b}}.
\end{array}
\end{equation}

\vskip 2mm

Since, after the first iteration of Kaczmarz's method (which ends with ${\mathbf{p}}_N$), we set
${\mathbf{x}}_1={\mathbf{p}}_N$, and we start the second iteration
from ${\mathbf{p}}_0={\mathbf{x}}_1$, it means that we are
continuing the original steps \eqref{kac}, and
that the new vector ${\mathbf{p}}_1$ is, in fact, the vector
${\mathbf{p}}_{N+1}$ of \eqref{kac}, that the new vector
${\mathbf{p}}_2$ is ${\mathbf{p}}_{N+2}$, and so on. Thus, Kaczmarz's method is only a renumbering of the single steps, as previously explained, keeping only those whose index is a multiple of
$N$, that is ${\mathbf{x}}_{n+1}={\mathbf{p}}_{(n+1)N}$ and
${\mathbf{r}}_{n+1}={\boldsymbol\varrho}_{(n+1)N}$. The main
interest of this renumbering lies in the relations \eqref{Qc}
which express the connection between two consecutive iterations by means of the fixed matrix $Q$, instead of
relations using matrices changing at each step of an iteration.
We similarly have
${\boldsymbol\varrho}_{(n+1)N+j-1}=Q^{(j)}{\boldsymbol\varrho}_{nN+j-1}$,
$j=1,\ldots,N$, with $Q^{(j)}=P_{j-1}\cdots P_1P_n\cdots P_j$.
Notice that $Q^{(1)}$ is identical to $Q$.

The matrices $P$ and $Q$ are similar, and it holds $\rho(Q) < 1$ as proved in \cite{gast}. Thus, the sequence
$({\mathbf x}_n)$ generated by Kaczmarz's method converges to the solution $\mathbf x$ of the system. Moreover, it follows from standard results on iterations of
the form \eqref{pqk} that
$$\|{\mathbf{x}}-{\mathbf{x}}_n\|={\cal O}(\rho(Q)^n).$$

By slightly improving Thm. 4.4 of \cite{gala}, which is based on a result by Meany \cite{mean} on the norm of a product of orthogonal projections of rank $N-1$, we have the

\begin{propo}
$$\|{\mathbf{x}}-{\mathbf{x}}_{n+1}\|^2 \leq \left( 1-\frac{(\det A)^2}{\prod_{i=1}^N \|A^T{\mathbf{e}}_i\|^2}\right) \|{\mathbf{x}}-{\mathbf{x}}_n\|^2.$$
\end{propo}

\vskip 2mm

A generalization of this result was recently given in \cite{ZZB} for the case
where $A$ and ${\mathbf{b}}$ are partitioned into blocks of rows.

Let us give some details about the block version.
Assume that the matrix $A$ is partitioned into the blocks of rows $A_1^T,A_2^T,\ldots$, where the matrix $A_i^T \in \rset^{N_i \times N}$ contains the rows $N_1+\cdots+N_{i-1}+1$ up to $N_1+\cdots+N_{i-1}+N_i$ of $A$  (with $N_0=0$), and the vector ${\mathbf{b}}$ is partitioned accordingly.
Each single step of the method now consists in the treatment of a block as a whole, and one complete cycle of all the blocks in their natural ordering is called the (cyclic) {\it block Kaczmarz's method}. As for Kaczmarz's method, we are able to give a matrix interpretation of this extension which
remains valid with now
${\boldsymbol\alpha}_i=(A_i^T)^\dag=A_i(A_i^TA_i)^{-1}$ and $P_i=I-A{\boldsymbol\alpha}_iE_i^T$, where
$E_i \in \rset^{N \times N_i}$ is the matrix whose columns are the vectors
$e_{N_1+\cdots+N_{i-1}+1},\ldots, e_{N_1+\cdots+N_{i-1}+N_i}$of the canonical basis of $\rset^N$, and it holds
$Q_i=I-A_i(A_i^TA_i)^{-1}A_i^T$. With these notations, one step of the block Kaczmarz's method writes
$${\mathbf p}_i={\mathbf p}_{i-1}+(A_i^T)^\dag E_i^T ({\mathbf b}-A{\mathbf p}_{i-1}).$$
This relation clearly generalizes \eqref{mmm}.

\section{Convergence acceleration}
\label{sacc}

For accelerating the convergence of a sequence, it can be transformed into another one which, under some assumptions, converges faster to the same limit. The idea behind such a {\it sequence transformation} is to assume that the sequence to be accelerated behaves as a model sequence, or satisfies some property, depending on unknowns parameters.
The set of these sequences is called the {\it kernel} of the transformation.
The unknown parameters are determined by imposing that the sequence interpolates the model sequence from a certain starting index $n$. Then, the limit of the model sequence, which depends on $n$ since the parameters depend on it, is taken as an approximation of the limit of the sequence to be accelerated which is thus transformed into a new sequence.
Of course, by construction, if the initial sequence belongs to the kernel of the transformation then, for all $n$, its limit is obtained. Although this observation was never proved rigourously, if the sequence is {\it close} in some sense to the kernel of the transformation, there is a good chance that it will be accelerated. This is why the notion of kernel is so important.

In practice, when having to accelerate a given sequence, one can use a known sequence transformation (also called an
{\it acceleration algorithm} or an {\it extrapolation method}), and verify that it can be accelerated by it.
Another approach is, starting from some algebraic property of the sequence to be accelerated, to construct a special transformation adapted to it. In both cases, the behavior of the sequence has to be analyzed or has to be characterized by some property.

On convergence acceleration methods for vector sequences, see, for instance, \cite{cbmrz}.

\subsection{The sequence generated by Kaczmarz's method}

Consider the sequence of vectors $({\mathbf{x}}_n)$ obtained by
Kaczmarz's method. Let $\Pi_\nu$ be the minimal
polynomial of the matrix $Q$ for the vector
${\mathbf{x}}-{\mathbf{x}}_0$, that is the polynomial of smallest
degree $\nu \leq N$ such that
$\Pi_\nu(Q)({\mathbf{x}}-{\mathbf{x}}_0)=0$. Since
$\Pi_\nu(Q)({\mathbf{x}}-{\mathbf{x}}_0)=A^{-1}\Pi_\nu(P)A({\mathbf{x}}-{\mathbf{x}}_0)=A^{-1}\Pi_\nu(P) \mathbf{r}_0$,
$\Pi_\nu$ is also the minimal polynomial of $P$ for the vector
$\mathbf{r}_0$. Setting $\Pi_\nu(\xi)=c_0+c_1\xi+\cdots+c_\nu\xi^\nu$, it
holds from \eqref{pqk}
\begin{equation}
Q^n\Pi_\nu(Q)({\mathbf{x}}-{\mathbf{x}}_0)=c_0({\mathbf{x}}-{\mathbf{x}}_n)+c_1({\mathbf{x}}-{\mathbf{x}}_{n+1})+
\cdots+c_\nu({\mathbf{x}}-{\mathbf{x}}_{n+\nu})=0,
\quad n=0,1,\ldots \label{ker}
\end{equation}
Thus, it follows that the vectors ${\mathbf{x}}-{\mathbf{x}}_n$
produced by Kaczmarz's method satisfy such an
homogeneous linear difference equation of order $\nu$, whose
solution is well--known.

If $Q$ is nondefective and if we denote by $\tau_1,\ldots,\tau_\nu$ the distinct zeros of $\Pi_\nu$,
this solution writes
$${\mathbf{x}}-{\mathbf{x}}_n=\sum_{i=1}^\nu d_i \tau_i^n {\mathbf{v}}_i, \quad n=0,1,\ldots,$$
where the $d_i$'s are scalars and the ${\mathbf{v}}_i$'s the
eigenvectors corresponding to the $\tau_i$'s. If $Q$ is defective,
the expression for ${\mathbf{x}}-{\mathbf{x}}_n$ still involves
powers of the eigenvalues, but it is more complicated.

\subsection{Sequence transformations}

The vector $\mathbf x$ can be exactly computed if we are able to build a sequence transformation whose kernel consists of sequences of the form \eqref{ker}.
However, since, in practical applications, $\nu$ could be quite large, we will restrict ourselves to building a sequence transformation whose kernel contains all sequences of the form
\begin{equation}
a_0^{(n)}({\mathbf{x}}-{\mathbf{x}}_n)+a_1^{(n)}({\mathbf{x}}-{\mathbf{x}}_{n+1})+\cdots+
a_k^{(n)}({\mathbf{x}}-{\mathbf{x}}_{n+k})=0, \quad n=0,1,\ldots, \label{eqdif}
\end{equation}
with $ k\leq \nu$. Thus, solving this equation for the vector $\mathbf x$, gives us an approximation of it denoted ${\mathbf{y}}_k^{(n)}$ since it depends on $k$ and $n$.

For that purpose, one has to find a procedure for computing the
unknown coefficients $a_0^{(n)},\ldots,a_k^{(n)}$, taking into account that they are
numbers while ${\mathbf{x}}-{\mathbf{x}}_{n+i}$ are vectors. All
such transformations can be considered as vector generalizations
of {\it Shanks' transformation} \cite{shanks} for sequences of
numbers. The common idea behind these transformations is to obtain
a system of linear equations whose solution is $a_0^{(n)},\ldots,a_k^{(n)}$, and then to compute ${\mathbf{y}}_k^{(n)}\simeq \mathbf x$.

\vskip 2mm

It turns out that several vector sequence transformations based on (or including) such a kernel already exist and have been studied
by various authors: the various $\varepsilon$--algorithms \cite{wynn,evec,etopo}, the Minimal
Polynomial Extrapolation (MPE) \cite{mpe}, the Modified Minimal  Polynomial Extrapolation (MMPE)
\cite{etopo,mmpe,sidi1}, the Reduced Rank Extrapolation (RRE) \cite{rre,rre1}, Germain--Bonne
transformations \cite{bgb}, the $H$--algorithm \cite{Halg}, and the $E$--algorithm \cite{ealg}.
They are described and analyzed, for example, in \cite{cbvec,etopo,rapdet,ealg,sidi2,sidi,sibr,smi,sidi1}.

\vskip 2mm

Writing \eqref{eqdif} for the indexes $n$ and $n+1$, subtracting, and multiplying scalarly by a vector $\mathbf y$ leads to the scalar equation
\begin{equation}
a_0^{(n)} ({\mathbf{y}},\Delta {\mathbf{x}}_n)+a_1^{(n)} ({\mathbf{y}},\Delta
{\mathbf{x}}_{n+1})+\cdots+a_k^{(n)}({\mathbf{y}},\Delta
{\mathbf{x}}_{n+k})=0, \label{yeq}
\end{equation}
where $\Delta$ is the difference operator defined by $\Delta {\mathbf{x}}_n ={\mathbf{x}}_{n+1} - {\mathbf{x}}_n.$

As explained in \cite[p. 39]{proj}, there are several possible strategies for constructing our system which,
apart from an additional normalization condition, has to contain $k$ equations
\begin{itemize}
\item use only one vector ${\mathbf{y}}$, and write \eqref{yeq}
for the indexes $n,\ldots,n+k-1$, \item write \eqref{yeq} only for
the index $n$, and choose $k$ linearly independent vectors
${\mathbf{y}}$,
\item write several relations \eqref{yeq}, and
choose several linearly independent vectors ${\mathbf{y}}$.
\end{itemize}

Adding to these $k$ equations the condition $a_0^{(n)}+\cdots+a_k^{(n)}=1$,
which does not restrict the generality (and is needed since the sum of the coefficients must be nonzero in order for
$\mathbf x$ in \eqref{eqdif} to be uniquely defined), we obtain a system of
$k+1$ equations in the $k+1$ unknowns $a_0^{(n)},\ldots,a_k^{(n)}$ (which
depend on $n$ and $k$). Its coefficients are denoted $d_{i,j}^{(n)}$,
and, for any of the preceding strategies, this system writes
\begin{equation}
\label{sys1}
\left\{ \begin{array}{l}
a_0^{(n)}+\cdots+a_k^{(n)}=1\\
d_{i,0}^{(n)}a_0^{(n)}+\cdots+d_{i,k}^{(n)}a_k^{(n)}=0,\qquad i=1,\ldots,k,
\end{array}
\right.
\end{equation}
where the coefficients $d_{i,j}^{(n)}$, for $i=1,\ldots,k$ and
$j=0,\ldots,k$, will be given below according to the chosen strategy. Then, for a fixed value of $k$, our sequence transformation
$({\mathbf{x}}_n) \longmapsto ({\mathbf{y}}_k^{(n)})$ is defined
by
\begin{equation}
{\mathbf{y}}_k^{(n)}=a_0^{(n)}{\mathbf{x}}_n+\cdots+a_k^{(n)}{\mathbf{x}}_{n+k},
\quad n=0,1,\ldots, \label{ykm1}
\end{equation}
which can be written as
\begin{equation}
\label{rapv}
{\mathbf{y}}_k^{(n)}=\frac{\left| \begin{array}{ccc}
           {\mathbf{x}}_{n} & \cdots & {\mathbf{x}}_{n+k} \\
           d_{1,0}^{(n)}& \cdots & d_{1,k}^{(n)}\\
           \vdots && \vdots\\
           d_{k,0}^{(n)}& \cdots & d_{k,k}^{(n)}
           \end{array}
                                  \right|}
                {\left| \begin{array}{ccc}
           1 & \cdots & 1 \\
d_{1,0}^{(n)}& \cdots & d_{1,k}^{(n)}\\
           \vdots && \vdots\\
           d_{k,0}^{(n)}& \cdots & d_{k,k}^{(n)}
           \end{array}
            \right|}=
          \frac{\left| \begin{array}{cccc}
           {\mathbf{x}}_{n} & \Delta {\mathbf{x}}_n &\cdots & \Delta {\mathbf{x}}_{n+k-1} \\
           d_{1,0}^{(n)}& \delta d_{1,0}^{(n)}& \cdots & \delta d_{1,k-1}^{(n)}\\
           \vdots && \vdots\\
           d_{k,0}^{(n)}& \delta d_{k,0}^{(n)}&\cdots & \delta d_{k,k-1}^{(n)}
           \end{array}
                                  \right|}
                {\left| \begin{array}{ccc}
\delta d_{1,0}^{(n)}& \cdots & \delta d_{1,k-1}^{(n)}\\
           \vdots && \vdots\\
           \delta d_{k,0}^{(n)}& \cdots & \delta d_{k,k-1}^{(n)}
           \end{array}
            \right|}
           \end{equation}
           where $\delta$ is the difference operator defined by
$\delta d_{i,j}^{(n)}=d_{i,j+1}^{(n)}-d_{i,j}^{(n)}$. The second determinantal expression shows that the vector ${\mathbf{y}}_k^{(n)}$ is the Schur complement \cite{proj}
$${\mathbf{y}}_k^{(n)}={\mathbf{x}}_n-[\Delta {\mathbf{x}}_n,\ldots,\Delta {\mathbf{x}}_{n+k-1}]\left(
\begin{array}{ccc}
\delta d_{1,0}^{(n)}& \cdots & \delta d_{1,k-1}^{(n)}\\
           \vdots && \vdots\\
           \delta d_{k,0}^{(n)}& \cdots & \delta d_{k,k-1}^{(n)}
           \end{array}
            \right)^{-1}
            \left(
            \begin{array}{c}
            d_{1,0}^{(n)}\\
            \vdots\\
            d_{k,0}^{(n)}
            \end{array}
            \right).$$
However, ${\mathbf{y}}_k^{(n)}$ is not a projection.

\vskip 2mm

The second determinantal formula in \eqref{rapv} means that all our sequences transformations can be also
written as
\begin{equation}
{\mathbf{y}}_k^{(n)}={\mathbf{x}}_n-\alpha_1^{(n)} \Delta {\mathbf{x}}_n
-\cdots-\alpha_k^{(n)} \Delta {\mathbf{x}}_{n+k-1}, \quad k=0,1,\ldots,
\label{ykm2}
\end{equation}
where the $\alpha_i^{(n)}$'s are solution of the linear system
\begin{equation}
\label{sys2}
\delta d_{i,0}^{(n)}\alpha_1^{(n)}+\cdots+\delta d_{i,k-1}^{(n)}\alpha_k^{(n)}=d_{i,0}^{(n)}, \quad i=1,\ldots,k.
\end{equation}

\vskip 2mm

Let us now specify various choices of the coefficients $d_{i,j}^{(n)}$,
for $i=1,\ldots,k$ and $j=0,\ldots,k$.

\subsubsection{The vector Shanks' transformations}
\label{vst}

We first consider sequence transformations which are directly inspired by the scalar sequence transformation of Shanks \cite{shanks}, and can be recursively implemented by Wynn's $\varepsilon$--algorithm \cite{wynn} or its generalizations, or by the $E$--algorithm \cite{ealg}, or the $H$--algorithm \cite{Halg}.

\begin{itemize}
\item Choosing $d_{i,j}^{(n)}=({\mathbf{e}}_p,\Delta
{\mathbf{x}}_{n+i+j-1})$, and replacing ${\mathbf{x}}_{j}$ in the
first row of the numerator of the first ratio in \eqref{rapv} by
$({\mathbf{e}}_p,{\mathbf{x}}_{j})$ corresponds to Shanks'
transformation applied componentwise to the vector sequence
$({\mathbf{x}}_n)$, and it gives the $p$th component of the
vector ${\mathbf{y}}_k^{(n)}$.
This transformation can be recursively implemented, separately for each component $p=1,\ldots,N$, by
the scalar $\varepsilon$--algorithm of Wynn \cite{wynn}.

\item The $\varepsilon$--algorithm can be also applied directly to a sequence of vectors after
defining the inverse ${\mathbf{u}}^{-1}$ of a vector
${\mathbf{u}}$ as
${\mathbf{u}}^{-1}={\mathbf{u}}/({\mathbf{u}},{\mathbf{u}})$. This
vector $\varepsilon$--algorithm was also proposed by Wynn \cite{evec}.
However, for this algorithm, \eqref{rapv} is no longer valid and the determinants have
to be replaced either by bigger and more complicated ones
\cite{prgm}, or by designants which generalize them in a
non--commutative algebra \cite{5sala,ahsa}. Thus, in this case, an
underlying system of linear equations for ${\mathbf{y}}_k^{(n)}$
does not exist, and this transformation has to be recursively
implemented by the vector $\varepsilon$--algorithm whose rules are, for a
sequence $(\mathbf u_n)$ of real (to simplify) vectors in $\rset^N$,
\begin{eqnarray*}
{\boldsymbol\varepsilon}^{(n)}_{-1} &=& 0 \in \rset^N, \quad
{\boldsymbol\varepsilon}^{(n)}_{0} = \mathbf{u}_{n} \in \rset^N, \quad
n=0,1, \ldots  \\
 {\boldsymbol\varepsilon}^{(n)}_{k+1} & = & {\boldsymbol\varepsilon}^{(n+1)}_{k-1}  +({\boldsymbol\varepsilon}^{(n+1)}_{k}
- {\boldsymbol\varepsilon}^{(n)}_{k})^{-1},
\quad k,n = 0,1, \ldots
\end{eqnarray*}

\item The choice $d_{i,j}^{(n)}=({\mathbf{y}},\Delta
{\mathbf{x}}_{n+i+j-1})$, where ${\mathbf{y}}$ is any nonzero vector
so that the denominators of \eqref{rapv} differs from zero, leads
to the topological Shanks' transformation introduced in
\cite{etopo}. It can be implemented via the topological
$\varepsilon$--algorithm. In this algorithm, the inverses are
defined in a different way for an even or an odd lower index as
$$
\begin{array}{lcl}
({\boldsymbol\varepsilon}^{(n+1)}_{2k}  -
{\boldsymbol\varepsilon}^{(n)}_{2k})^{-1}&=&{\mathbf{y}}/({\mathbf{y}},{\boldsymbol\varepsilon}^{(n+1)}_{2k}
- {\boldsymbol\varepsilon}^{(n)}_{2k}) \\
({\boldsymbol\varepsilon}^{(n+1)}_{2k+1}  -
{\boldsymbol\varepsilon}^{(n)}_{2k+1})^{-1}&=& ({\boldsymbol\varepsilon}_{2k}^{(n+1)}-
{\boldsymbol\varepsilon}_{2k}^{(n)})/ ({\boldsymbol\varepsilon}^{(n+1)}_{2k+1}  -
{\boldsymbol\varepsilon}^{(n)}_{2k+1}, {\boldsymbol\varepsilon}^{(n+1)}_{2k}  -
{\boldsymbol\varepsilon}^{(n)}_{2k}).
\end{array}
$$


\end{itemize}

In the preceding $\varepsilon$-algorithms, only the vectors (or the numbers) with
an even lower index ${\boldsymbol\varepsilon}_{2k}^{(n)}$ are interesting for
the purpose of convergence acceleration, those with an odd lower
index ${\boldsymbol\varepsilon}_{2k+1}^{(n)}$ being intermediate computations.
Applying any of them to the sequence
$({\boldsymbol\varepsilon}^{(n)}_{0} ={\mathbf{x}}_n)$ produced by Kaczmarz's method gives
${\boldsymbol\varepsilon}_{2k}^{(n)}={\mathbf{y}}_k^{(n)}$.
Notice that the computation of one vector ${\boldsymbol\varepsilon}_{2k}^{(n)}$
requires the $2k+1$ vectors ${\mathbf{x}}_n,\ldots,{\mathbf{x}}_{n+2k}$. Thus, when $k$ is increased by 1, the computation of each new vector ${\boldsymbol\varepsilon}_{2k}^{(n)}$
needs two additional iterates of Kaczmarz's method, while it only requires one when $n$ increases.




\subsubsection{The MPE, MMPE, and RRE}
\label{mmr}

The following choices also lead to known transformations but, in this case,  the computation of ${\mathbf{y}}_k^{(n)}$ only requires the $k+2$ vectors ${\mathbf{x}}_n,\ldots,{\mathbf{x}}_{n+k+1}$.
Algorithms for their recursive implementation also exist \cite{bvt,Halg,fs}.

\begin{itemize}

\item The choice $d_{i,j}^{(n)}=({\mathbf y}_i, \Delta {\mathbf
x}_{n+j})$, for $i=1,\ldots,k$ and $j=0,\ldots,k$ corresponds to
the MMPE (Modified Minimal Polynomial Extrapolation)
\cite{mmpe,sidi1,etopo}, where the
${\mathbf{y}}_i$ are linearly independent vectors. This
transformation is another generalization of the topological Shanks
transformation also given in \cite{etopo}.

\item If we take $d_{i,j}^{(n)}=(\Delta {\mathbf{x}}_{n+i-1},\Delta
{\mathbf{x}}_{n+j})$, for $i=1,\ldots,k$ and $j=0,\ldots,k$, we
obtain the MPE (Minimal Polynomial Extrapolation) \cite{mpe}. This method is
mathematically equivalent to a transformation due to  Germain--Bonne
\cite{bgb}. Other choices are proposed in the same reference.

\item The choice $d_{i,j}^{(n)}=(\Delta^2 {\mathbf{x}}_{n+i-1},\Delta
{\mathbf{x}}_{n+j})$ leads to the RRE
(Reduced Rank Extrapolation) \cite{rre,rre1}.

\end{itemize}

It must be noticed that the vector ${\mathbf{y}}_k^{(n)}$ always exists for the RRE but not for the MPE or the MMPE.
Existence conditions are discussed in \cite{js} and \cite{sidi6}.

\subsubsection{The simplest transformations}

For $k=1$, the transformations have the following very simple forms

\begin{itemize}
\item the MMPE and the topological Shanks' transformation write
$${\mathbf{y}}_1^{(n)}={\mathbf{x}}_n-\frac{({\mathbf{y}}_1,\Delta {\mathbf{x}}_n)}{({\mathbf{y}}_1,\Delta^2 {\mathbf{x}}_n)}\Delta {\mathbf{x}}_n,$$
\item for the MPE and the transformation due to Germain--Bonne, we have
$${\mathbf{y}}_1^{(n)}={\mathbf{x}}_n-\frac{(\Delta {\mathbf{x}}_n,\Delta {\mathbf{x}}_n)}{(\Delta {\mathbf{x}}_n,\Delta^2 {\mathbf{x}}_n)}\Delta {\mathbf{x}}_n,$$
\item for the RRE, it holds
$${\mathbf{y}}_1^{(n)}={\mathbf{x}}_n-\frac{(\Delta^2 {\mathbf{x}}_n,\Delta {\mathbf{x}}_n)}{(\Delta^2 {\mathbf{x}}_n,\Delta^2 {\mathbf{x}}_n)}\Delta {\mathbf{x}}_n,$$
which is the same as the one proposed in \cite{gea},
\item for the vector $\varepsilon$--algorithm, we obtain
$${\mathbf{y}}_1^{(n)}={\mathbf{x}}_{n+1}+\frac{{\boldsymbol{\varepsilon}}_1^{(n+1)}-{\boldsymbol{\varepsilon}}_1^{(n)}}
{({\boldsymbol{\varepsilon}}_1^{(n+1)}-{\boldsymbol{\varepsilon}}_1^{(n)},
{\boldsymbol{\varepsilon}}_1^{(n+1)}-{\boldsymbol{\varepsilon}}_1^{(n)})} \mbox{~~with~~}
{\boldsymbol{\varepsilon}}_1^{(n)}=\frac{\Delta {\mathbf{x}}_{n}}
{(\Delta {\mathbf{x}}_{n},\Delta {\mathbf{x}}_{n})}.$$
\end{itemize}

\subsubsection{Orthogonality properties}

Using the determinantal formula \eqref{rapv}, we have the following orthogonality properties

\begin{propo}
~~\\
The vector
$Q({\mathbf{y}}_k^{(n)}-{\mathbf{x}})-({\mathbf{y}}_k^{(n)}-{\mathbf{x}})$
is orthogonal to \vskip 1mm

$\bullet$ $\Delta {\mathbf{x}}_n,\ldots,\Delta {\mathbf{x}}_{n+k-1}$ for the MPE (which is identical to
Germain--Bonne transformation),

$\bullet$ ${\mathbf{y}}$ for the topological Shanks transformation
(topological $\varepsilon$--algorithm),

$\bullet$ $\Delta^2 {\mathbf{x}}_n,\ldots,\Delta^2
{\mathbf{x}}_{n+k-1}$ for the RRE,

$\bullet$ ${\mathbf{y}}_1,\ldots,{\mathbf{y}}_k$ for the MMPE.
\end{propo}

For a discussion of the properties of these algorithms and their application to the solution of systems of linear and nonlinear equations, see \cite{js}.

\section{Acceleration of Kaczmarz's method}

Although they are different, all the acceleration methods
presented in Section \ref{sacc} have the same kernel, namely
sequences satisfying \eqref{eqdif}. Thus, when applied to any sequence of vectors
generated by iterations of the form
${\mathbf{x}}_{n+1}=B{\mathbf{x}}_n+{\mathbf{c}}$, $n=0,1,\ldots$,
we obtain ${\mathbf{y}}_{\nu}^{(n)}=(I-B)^{-1}{\mathbf{c}}$ for all
$n$, where $\nu$ is the degree of the minimal polynomial of $B$
for the vector ${\mathbf{x}}-{\mathbf{x}}_0$. This result means that
these methods are direct methods for solving systems of linear equations
\cite{etopo}. Therefore, this property holds for the
iterates of Kaczmarz's method, since $\Pi_\nu(1) \neq 0$
(otherwise the matrix $I-Q$ would be singular). Of course, in
practice, these algorithms cannot be used for obtaining the exact
solution when the dimension of the system is large since, for computing
${\mathbf{y}}_{\nu}^{(n)}$,
 they
require the storage of too many vectors.

\subsection{The AK and RK algorithms}

The preceding algorithms can be used in two different ways for
accelerating the iterations \eqref{cycle} produced by
Kaczmarz's method as it will now be explained. Since the
computation of ${\mathbf{y}}_k^{(n)}$ does not require the same
number of vectors issued from Kaczmarz's method,
according to the transformation used, in order to have a
unified presentation, we will denote by $\ell+1$ the number of
vectors ${\mathbf{x}}_i$ needed to compute ${\mathbf{y}}_k^{(n)}$
by any of the preceding procedures. Remember that $\ell$ depends
on $k$ and that $\ell=2k$ for the $\varepsilon$--algorithms, and
$\ell=k+1$ for the other algorithms.

\vskip 2mm

\begin{itemize}
\item The first one consists to apply one of the algorithms
implementing a sequence transformation to the sequence
${\mathbf{x}}_0,{\mathbf{x}}_1,\ldots$ given by
Kaczmarz's method, and, after fixing the index $k$, to build
simultaneously the sequence
${\mathbf{z}}_0={\mathbf{y}}_k^{(0)},{\mathbf{z}}_1={\mathbf{y}}_k^{(1)},\ldots$.
The computation of ${\mathbf{y}}_k^{(0)}$ can only begin after having computed the iterate
${\mathbf{x}}_\ell$, but the computation of each new
transformed vector needs only one new iterate of Kaczmarz's method.
This procedure is called the {\it accelerated Kaczmarz} (AK)
algorithm. Let us give the general structure for the implementation of the
AK algorithm.

\begin{center}
\hrulefill \\[0pt]
\vspace{-0.1cm} \textbf{Accelerated Kaczmarz (AK) algorithm} \\[0pt]
\vspace{-0.2cm} \hrulefill \\[0pt]
\begin{tabbing}
NN \= XX \= XX \= XX \= XX \= \kill {\bf Require} $A\in
{\mathbb{R}}^{N\times N},\;{\mathbf{b}}\in {\mathbb{R}}^{N}$,
${\mathbf{x}}_0\in {\mathbb{R}}^{N}$\\
{\bf Choose} $k \in {\mathbb{N}}, k \geq 1$\\
{\bf Set} $\ell = k+1$ or $\ell = 2k$\\
{\bf for} $n=0,1,\ldots$ {\rm until convergence} {\bf do}\\
\> ${\mathbf{p}}_{0} \leftarrow {\mathbf{x}}_{n} $\\
\> {\bf Compute} ${\mathbf{p}}_{i}, \quad i=1, \ldots,N$\\
\> ${\mathbf{x}}_{n+1} \leftarrow {\mathbf{p}}_{N} $\\
\> {\bf If} $n \geq \ell-1$ {\bf then}\\
\> \> {\bf Compute} $\mathbf{y}_{k}^{(n-\ell+1)}$ \\
\> \> ${\mathbf{z}}_{n-\ell+1} \leftarrow \mathbf{y}_{k}^{(n-\ell+1)}$\\
\> {\bf end if} \\
{\bf end for} $n$\\
\end{tabbing}
\vspace{-0.4cm} \hrulefill
\end{center}

\item In the second way, we set ${\mathbf{x}}_0$, and we compute
${\mathbf{x}}_1,\ldots,{\mathbf{x}}_{\ell}$ by
Kaczmarz's method, we apply one of the algorithms implementing a
sequence transformation to them, and we obtain
${\mathbf{z}}_0={\mathbf{y}}_k^{(0)}$. Then, we restart
Kaczmarz's method from ${\mathbf{z}}_0$, that is we set
${\mathbf{x}}_0={\mathbf{z}}_0={\mathbf{y}}_k^{(0)}$, we compute
the new $\ell$ vectors ${\mathbf{x}}_1,\ldots,{\mathbf{x}}_{\ell}$
by Kaczmarz's method, we apply again the acceleration algorithm to
these vectors ${\mathbf{x}}_0,\ldots,{\mathbf{x}}_{\ell}$, we
obtain ${\mathbf{z}}_1={\mathbf{y}}_k^{(0)}$, we restart
Kaczmarz's method from ${\mathbf{x}}_0={\mathbf{z}}_1$, and so on.
This second procedure is called the {\it restarted Kaczmarz} (RK)
algorithm and it can be implemented as follows.

\begin{center}
\hrulefill \\[0pt]
\vspace{-0.1cm} \textbf{Restarted Kaczmarz (RK) algorithm} \\[0pt]
\vspace{-0.2cm} \hrulefill \\[0pt]
\begin{tabbing}
NN \= XX \= XX \= XX \= XX \= \kill
{\bf Require} $A\in {\mathbb{R}}^{N\times N},\;{\mathbf{b}}\in {\mathbb{R}}^{N}$, ${\mathbf{x}}_0 \in {\mathbb{R}}^{N}$\\
{\bf Choose} $k \in {\mathbb{N}}, k \geq 1$\\
{\bf Set} $\ell = k+1$ or $\ell = 2k$\\
{\bf for} $n=0,1,\ldots$ {\rm until convergence} {\bf do}\\
\>{\bf for} $j=0,\ldots,\ell-1$  {\bf do}\\
\> \> ${\mathbf{p}}_{0} \leftarrow {\mathbf{x}}_{j} $\\
\> \> {\bf Compute} ${\mathbf{p}}_{i}, \quad i=1, \ldots,N$\\
\> \> ${\mathbf{x}}_{j+1} \leftarrow {\mathbf{p}}_{N} $\\
\> {\bf end for} $j$\\
\> {\bf Compute} $\mathbf{y}_k^{(0)}$ \\
\>  ${\mathbf{z}}_{n} \leftarrow \mathbf{y}_k^{(0)}$\\
\>  ${\mathbf{x}}_{0} \leftarrow \mathbf{z}_{n}$\\
{\bf end for} $n$
\end{tabbing}
\vspace{-0.4cm} \hrulefill
\end{center}

\end{itemize}

\subsection{The case of the $\varepsilon$--algorithms}
\label{qqq}

Let us now consider the particular case of the $\varepsilon$--algorithms.

If we assume that the eigenvalues of the matrix $Q$, defined in
Section \ref{MatInt}, are numbered such that
$|\tau_1|>|\tau_2|>\cdots>|\tau_N|>0$, then the sequences
$({\boldsymbol\varepsilon}_{2k}^{(n)})$, for $k \leq \nu$ fixed,
constructed by these algorithms are such that \cite{sidi,cbv,smi}
$$\|{\mathbf{x}}-{\boldsymbol\varepsilon}_{2k}^{(n)}\|={\cal O}(\tau_{k+1}^n), \quad n \to \infty.$$
If the preceding assumptions on the $\tau_i$'s are not satisfied
(which corresponds to several eigenvalues of $Q$ having the same
modulus or to a  defective matrix $Q$), the expression for
${\mathbf{x}}-{\mathbf{x}}_n$ is more complicated than given above
\cite{sibr} (see \cite{cbmc} for the complete expression in the
scalar case), but the convergence is still accelerated by the
$\varepsilon$--algorithms. These results also hold for the
topological $\varepsilon$--algorithm, independently of the choice
of the arbitrary vector ${\mathbf{y}}$ occurring in this
algorithm. Thus, each sequence $({\boldsymbol\varepsilon}_{2k}^{(n)})$ obtained
by the accelerated Kaczmarz algorithm converges to
${\mathbf{x}}$ faster than the preceding sequence
$({\boldsymbol\varepsilon}_{2k-2}^{(n)})$ when $n$ tends to infinity.
Quite similar results hold for the RRE and the MPE \cite{sidi6}.

\vskip 2mm

If the topological $\varepsilon$--algorithm is applied to a
sequence $({\mathbf{x}}_n)$ obtained by any iterative method of
the form ${\mathbf{x}}_{n+1}=B{\mathbf{x}}_n+{\mathbf{c}}$, with
the choice
${\mathbf{y}}={\mathbf{r}}_0={\mathbf{b}}-A{\mathbf{x}}_0$, then
the vectors ${\boldsymbol\varepsilon}_{2k}^{(0)}$, $k=0,1,\ldots$, are
identical to those obtained by Lanczos' method (that is, for
example, by the biconjugate algorithm or the conjugate gradient
algorithm in the symmetric positive definite case) \cite[Thms. 4.1
and 4.2, p. 186--7]{birk}. This property
comes out from the determinantal expressions of the vectors produced by the topological $\varepsilon$--algorithm and by
Lanczos' method. From \eqref{Qcc}, we
immediately see that this is the case of the vectors given by
Kaczmarz's method which correspond to $B=Q$ and
${\mathbf{c}}=A^{-1}(I-P){\mathbf{b}}$. Thus, if the topological
$\varepsilon$--algorithm is applied, with ${\mathbf{y}}={\mathbf{r}}_0$, to the sequence $({\mathbf{x}}_n)$ given by Kaczmarz's method, then the sequence
$({\boldsymbol\varepsilon}_{2k}^{(0)})$ is exactly the sequence produced by
Lanczos' method. Using it as explained in the restarted
Kaczmarz algorithm corresponds to restarting Kaczmarz's
method with the vector obtained by the Lanczos' method after $k$
of its iterates. At each restart, the vector ${\mathbf{y}}$ has to
be taken equal to the corresponding residual.

\section{Implementation}

The implementation of our convergence acceleration procedures can be realized by three different ways.
The first one consists, for a fixed value of $k$, in solving the linear system \eqref{sys1}
and using \eqref{ykm1}. The second way is to employ \eqref{ykm2} and the system \eqref{sys2}. The third possibility is to apply, when it exists,  a recursive algorithm. Since, for the vector
$\varepsilon$--algorithm, no underlying system of linear equations
is known, its implementation can only be realized through its
recursive rules, given in Section \ref{vst}.

\vskip 2mm

An important practical problem is the choice of $k$. On one side,
the effectiveness of our procedures seems to increase with $k$ but, on the
other side, the number of vectors to be stored also increases with it.
Thus,  if the system is
large, the value of $k$ has to be kept quite small since too
many vectors will have to be stored.
The numerical stability of the procedures has also to be
taken into consideration when $k$ increases.

Since
Kaczmarz's method often converges slowly, the quantities $\Delta
{\mathbf{x}}_{n+i}$ are small, and the elements of the linear systems
to be solved for obtaining the coefficients of the linear
combination giving ${\mathbf{y}}_k^{(n)}$ approach zero. Thus,
whatever the method used, the linear systems \eqref{sys1} or \eqref{sys2} are ill--conditioned
even for small values of $k$, and rounding errors can degrade the
results.
We tried several ways for computing their solution but the influence on the vectors
${\mathbf{y}}_k^{(n)}$ was small.

Instead of solving the systems, it is possible to use a recursive algorithm for the computation of
the vectors
${\mathbf{y}}_k^{(n)}$.
Usually, sequence transformations are
more efficiently implemented via a recursive algorithm as those
mentioned above, and, in our case, this kind of approach seems
to give results less sensitive to rounding errors for the topological and the vector $\varepsilon$--algorithms.

The quantities computed by all the recursive algorithms for implementing our transformations
can be displayed in a triangular array (or a lozenge one for the $\varepsilon$--algorithms) , and
${\mathbf{y}}_k^{(n)}$ corresponds to the lowest rightmost element of the
array. Its computation imposes to compute first all the preceding
elements in the array. Thus, it could be much too expensive to
store all these vectors if $k$ is not quite small and if $N$ is
large. Hopefully, it is possible to proceed in the array by ascending diagonals (that is, starting from ${\mathbf{y}}_0^{(n+k)}={\mathbf x}_{n+k}$, to compute
${\mathbf{y}}_1^{(n+k-1)},\ldots,{\mathbf{y}}_k^{(n)}$), and, thus, storing only one ascending
diagonal. For more details on such a technique, see \cite{cbmrz}.

\vskip 2mm
Another  problem is the choice of the arbitrary vector
$\mathbf y$ in the topological Shanks transformation, and of the vectors
${\mathbf{y}}_1,\ldots,{\mathbf{y}}_k$ appearing in the MMPE (see \cite{js} for a discussion about this choice). This is an unsolved problem and further studies
will be necessary.

\vskip 2mm

In the literature, the system $A{\mathbf x}={\mathbf b}$ to be solved by Kaczmarz' method is
often replaced by the system $DA{\mathbf x}=D{\mathbf b}$, where $D=\mbox{diag}(1/\|{\mathbf a}_i\|)$.
Thus, each row vector of the new matrix $DA$ will have a norm equal to 1.
This preconditioning, of course, simplify the computation of
${\mathbf p}_i$ in \eqref{mmm}, and we will use it.

\vskip 2mm

As a last comment, notice that, although some transformations are
equivalent from the theoretical point of view (for instance the
MPE and Germain--Bonne transformation, or the MMPE and the
$H$--algorithm),  they do not produce exactly the same
numerical results.   The same is true when we compare the results
of a transformation implemented in the different ways
described in this paper.

\section{Numerical results}

Our numerical experiments were performed using Matlab$^\circledR$ 7.11 and matrices coming out from the
{\tt gallery} set.
The solution was set to $\mathbf x=(1,\ldots,1)^T$ and the right hand side $\mathbf b$ was computed accordingly.
We took random vectors with components uniformly distributed in $[-1,+1]$ for the arbitrary vector
${\mathbf{y}}$ in the topological
$\varepsilon$--algorithm, and for the vectors ${\mathbf{y}}_i$ in
the MMPE.

\vskip 2mm

Some figures show the Euclidean norms of the errors and,
in order to compare the acceleration brought by each procedure, we
also give the ratios of the norms of the errors between the
iterate ${\mathbf{z}}_n$ obtained by the AK or the RK
algorithm and the iterate of Kaczmarz' method with the highest
index used in its construction (for AK), or the iterate with the highest index  which would have been used if we have led the method continue without restarting it (for RK) , that is
$$\frac{\|{\mathbf{z}}_n-{\mathbf{x}}\|}{\|{\mathbf{x}}_{n+\ell}-{\mathbf{x}}\|}~~\mbox{(AK)} \quad \mbox{and} \quad
\frac{\|{\mathbf{z}}_n-{\mathbf{x}}\|}{\|{\mathbf{x}}_{(n+1)(\ell+1)}-{\mathbf{x}}\|}~~\mbox{(RK)}.$$

\begin{figure}[hbt]
\begin{tabular}{@{}c@{}c@{}}
 \includegraphics[width=0.5\textwidth]{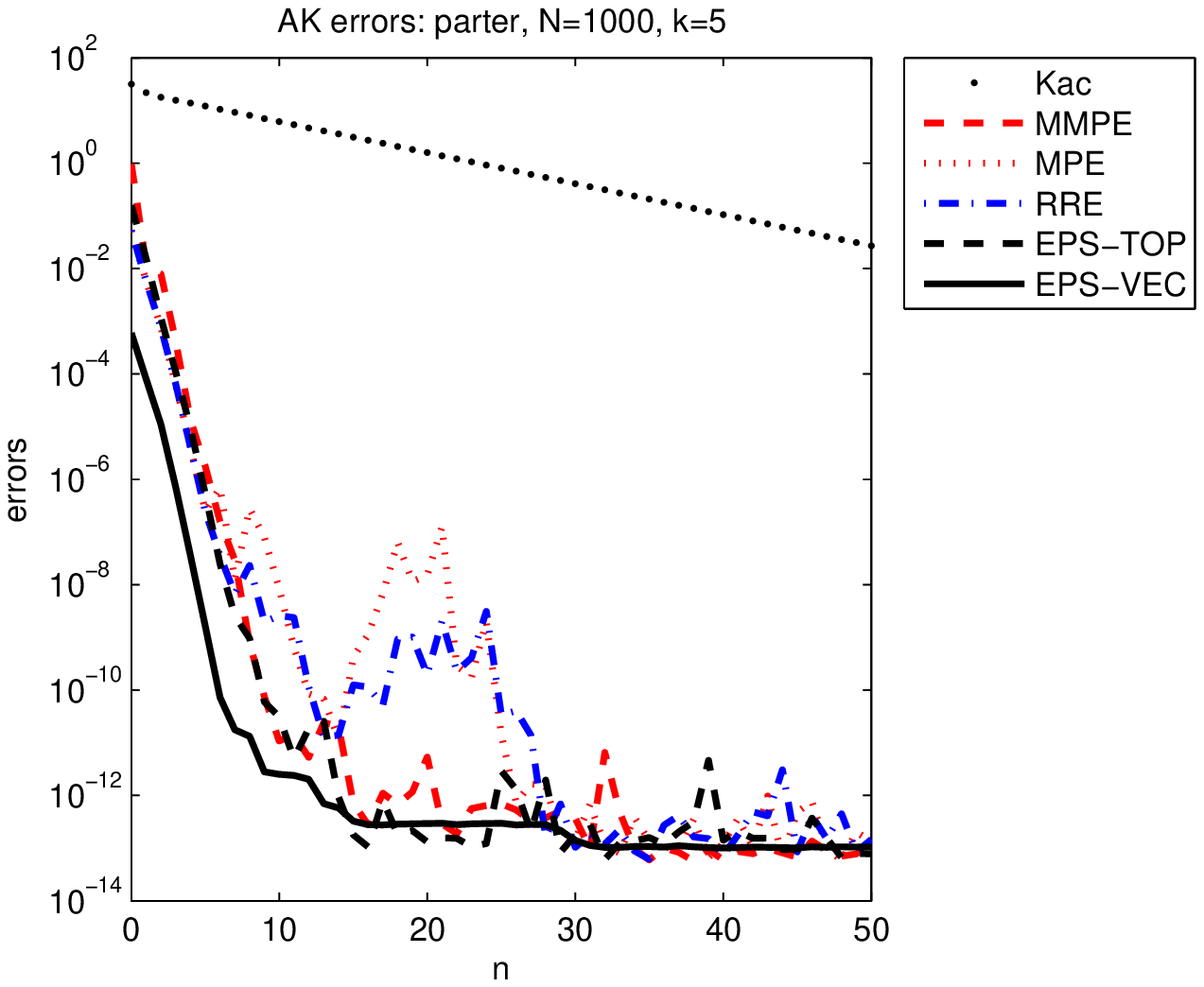}
 & \includegraphics[width=0.5\textwidth]{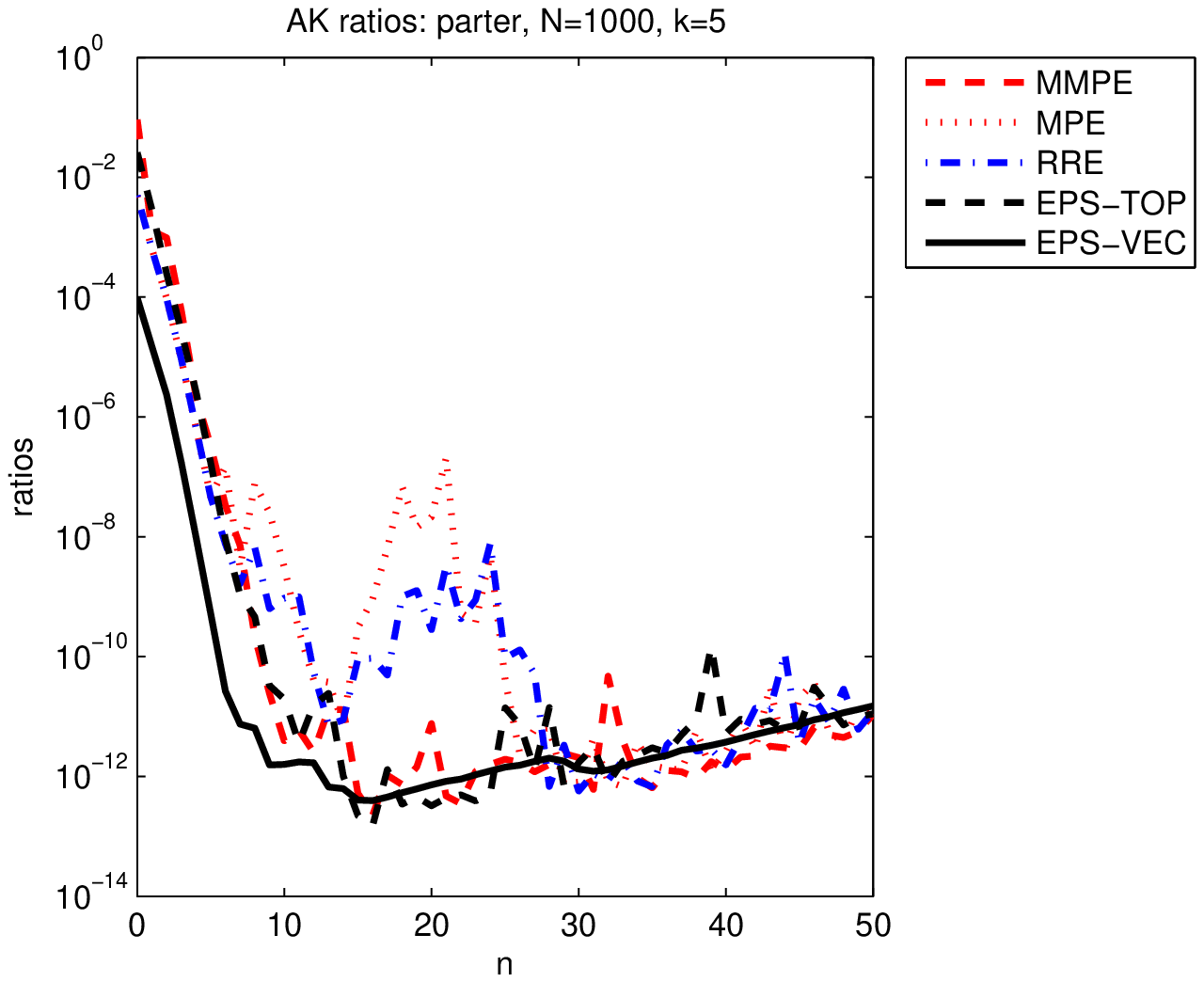}
 \end{tabular}
 \vspace{-0.7cm}
\caption{AK algorithm: errors and ratios for {\tt parter} matrix, $N=1000$, $k=5$.}
\label{parterAK}
\end{figure}

\begin{figure}[t!]
\begin{tabular}{@{}c@{}c@{}}
 \includegraphics[width=0.5\textwidth]{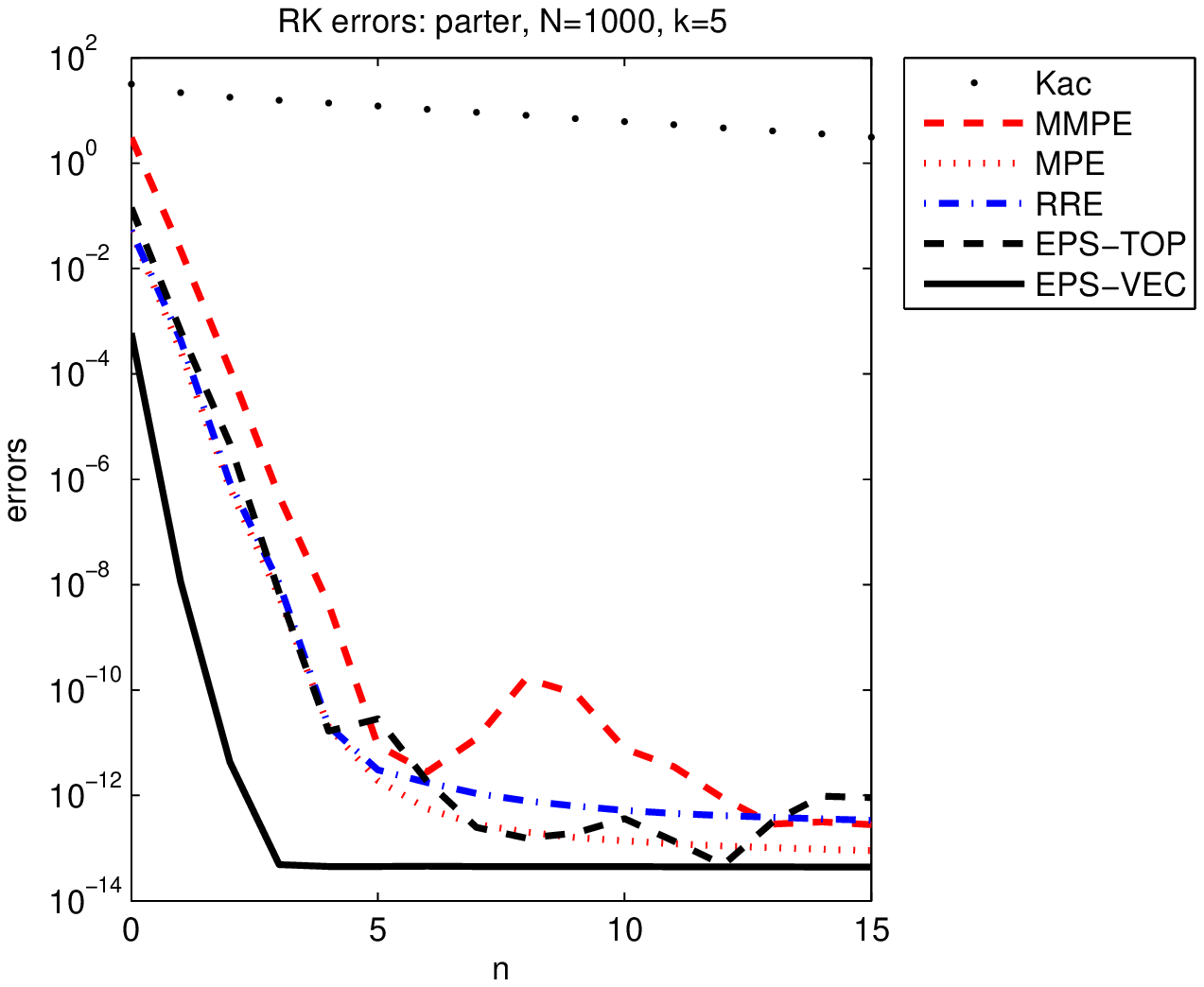}
 & \includegraphics[width=0.5\textwidth]{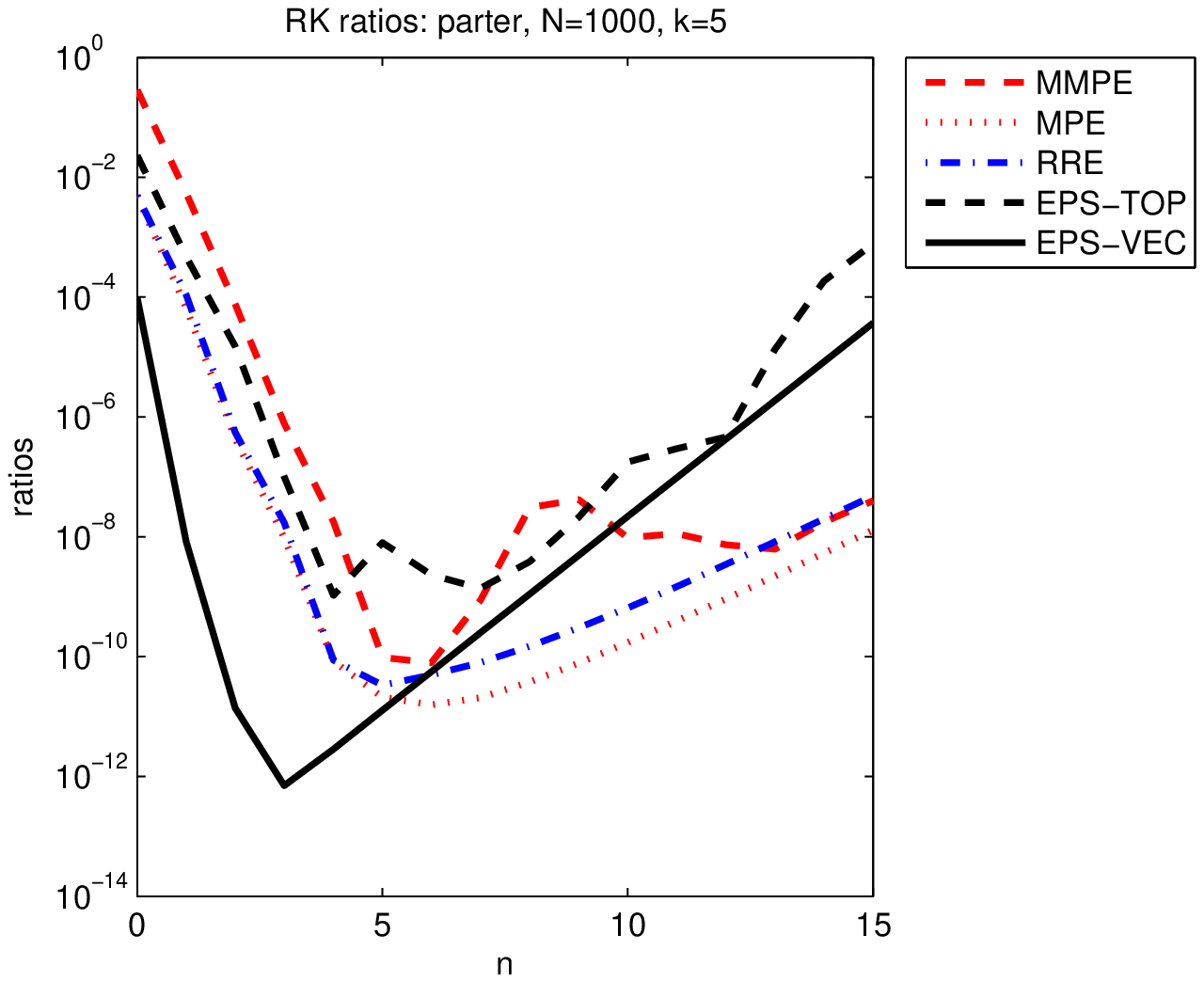}
 \end{tabular}
 \vspace{-0.7cm}
\caption{RK algorithm: errors and ratios for {\tt parter} matrix, $N=1000$, $k=5$.}
\label{parterRK}
\end{figure}

We consider the {\tt parter} matrix $A$, $N=1000$, $\kappa(A)\simeq 4.2306$, a Toeplitz matrix with singular values near $\pi$.
In the Figures \ref{parterAK} and \ref{parterRK}, we compare Kaczmarz' method with the
MMPE, MPE, RRE and the topological $\varepsilon$--algorithm, implemented by solving the system \eqref{sys1}, and
the
vector $\varepsilon$--algorithm, respectively for the AK and RK algorithms, with $k=5$.


\begin{figure}[t!]
\begin{tabular}{@{}c@{}c@{}}
 \includegraphics[width=0.5\textwidth]{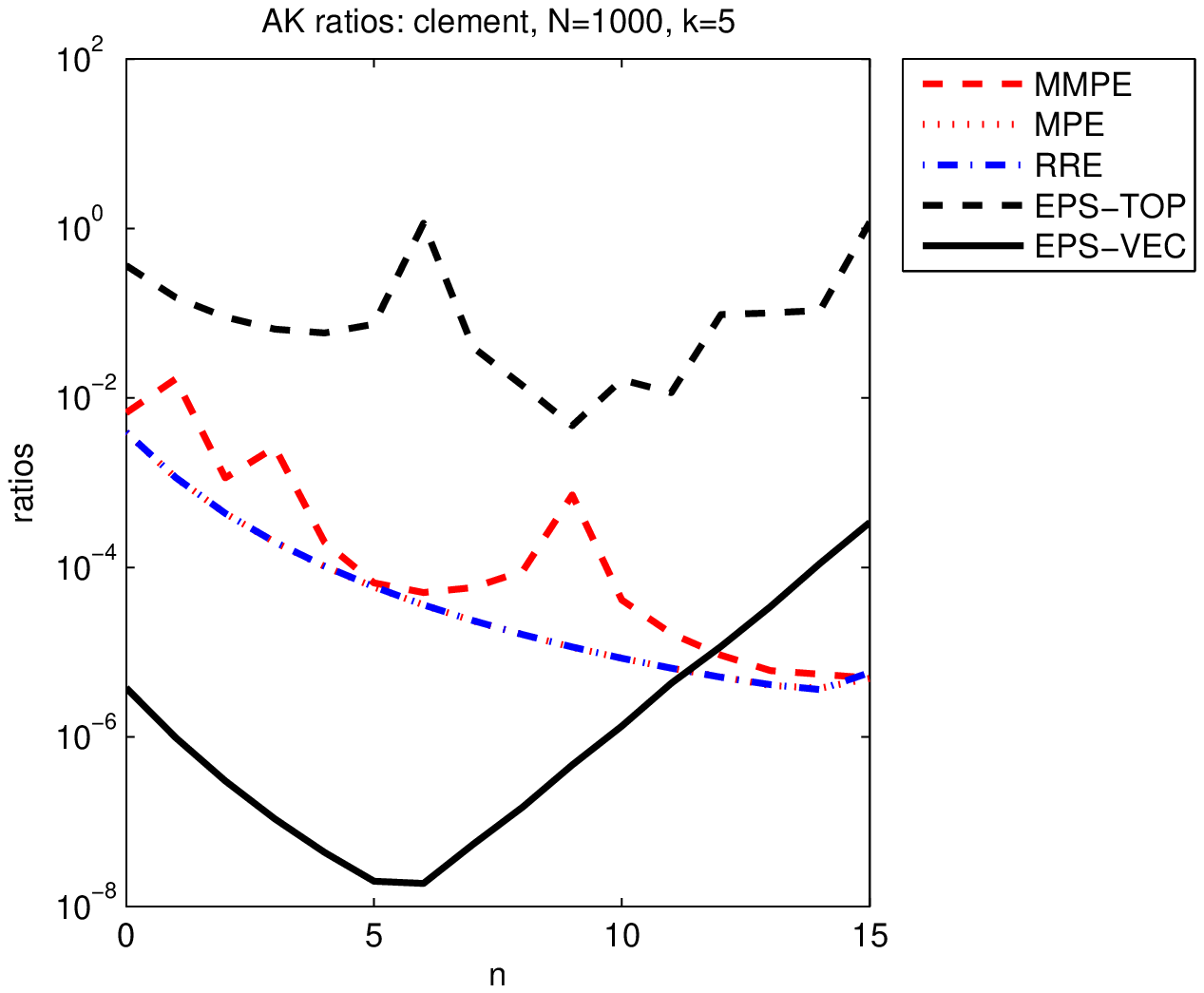}
 & \includegraphics[width=0.5\textwidth]{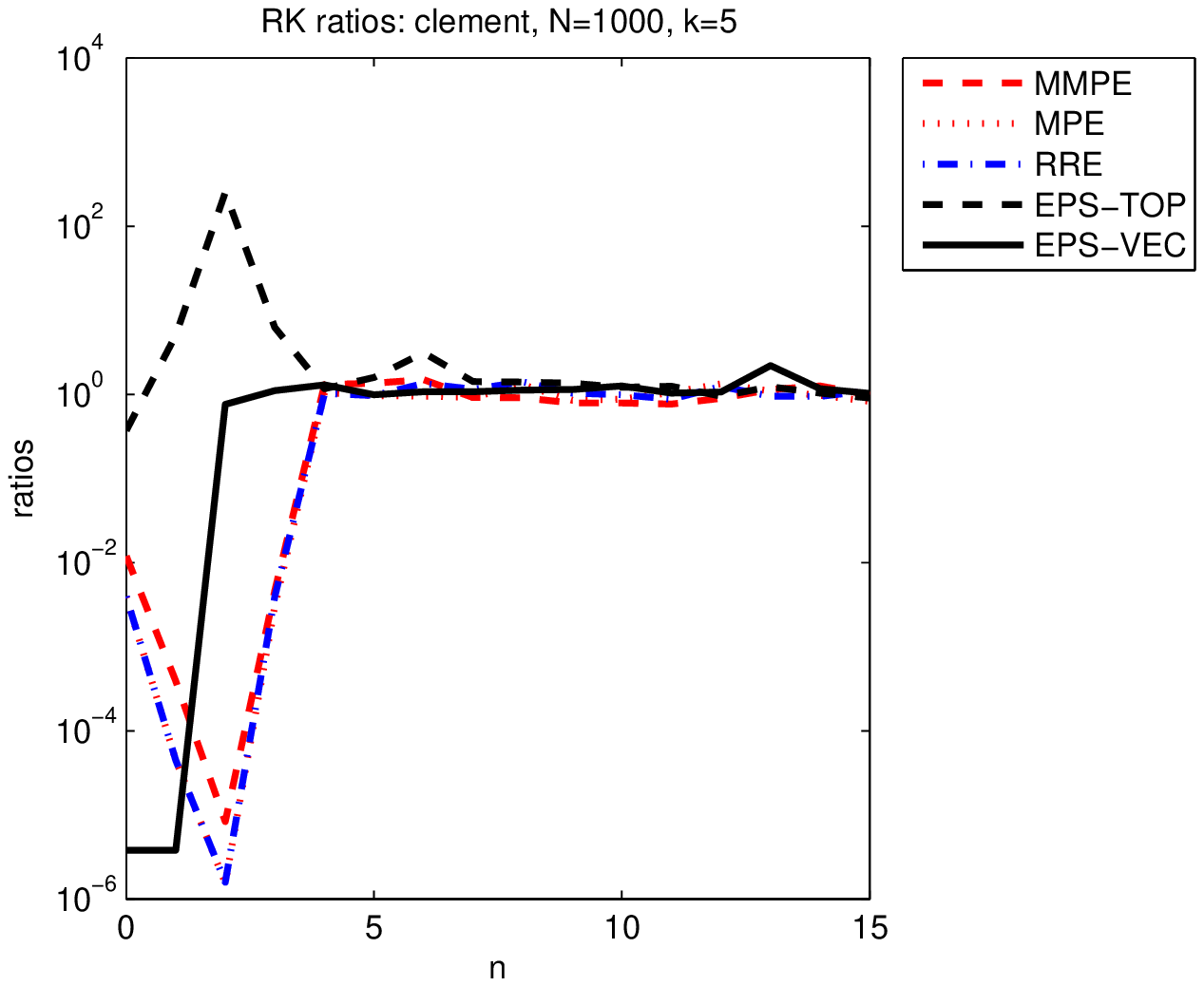}
 \end{tabular}
 \vspace{-0.7cm}
\caption{AK and RK algorithms: ratios for {\tt clement} matrix, $N=1000$, $k=5$.}
\label{clementAKRK}
\vspace{0.3cm}
\begin{tabular}{@{}c@{}c@{}}
 \includegraphics[width=0.5\textwidth]{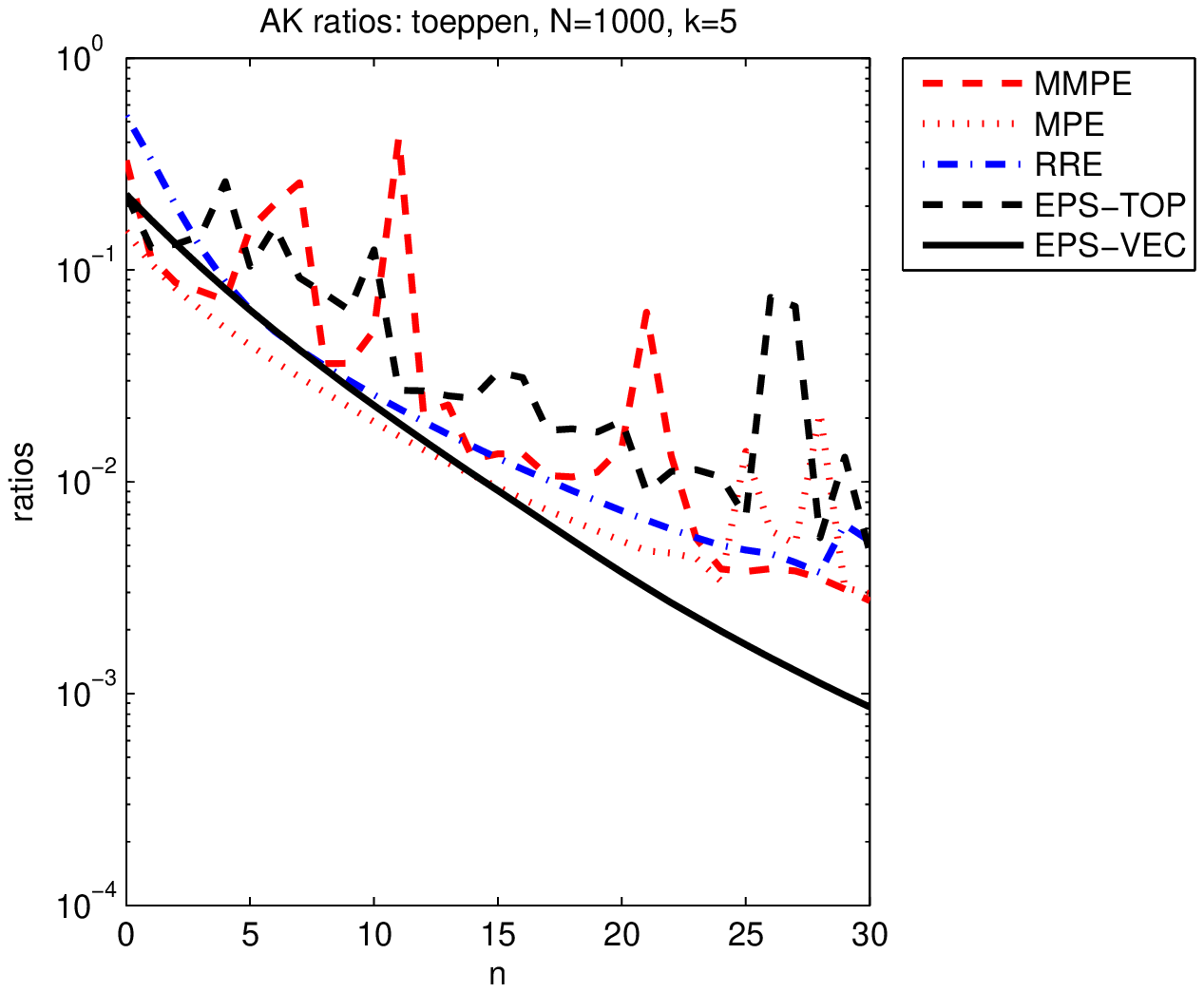}
 & \includegraphics[width=0.5\textwidth]{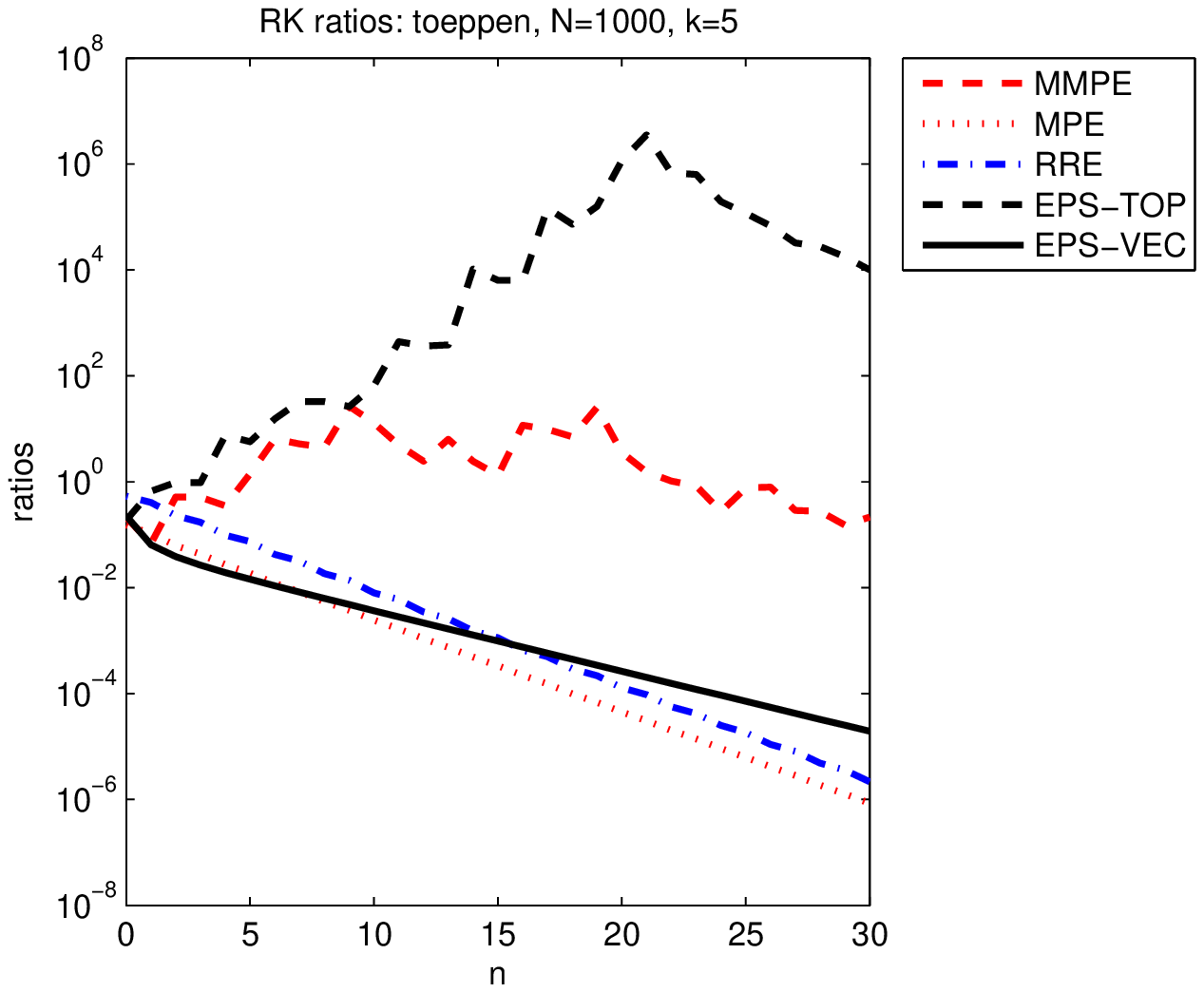}
 \end{tabular}
 \vspace{-0.7cm}
\caption{AK and RK algorithms: ratios for {\tt toeppen} matrix, $N=1000$, $k=5$.}
\label{toeppen}
\end{figure}

\begin{figure}[t!]
\begin{tabular}{@{}c@{}c@{}}
 \includegraphics[width=0.5\textwidth]{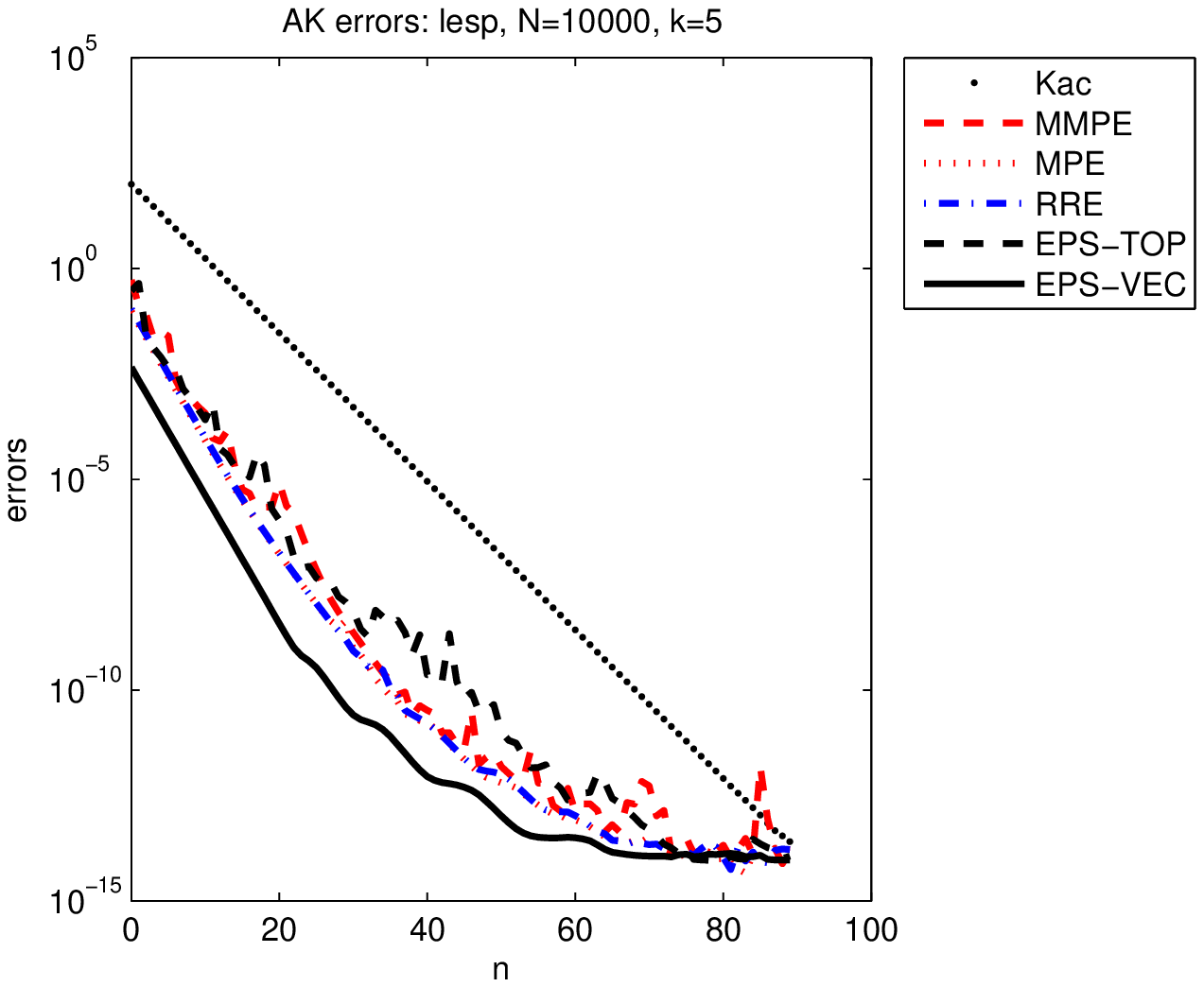}
 & \includegraphics[width=0.5\textwidth]{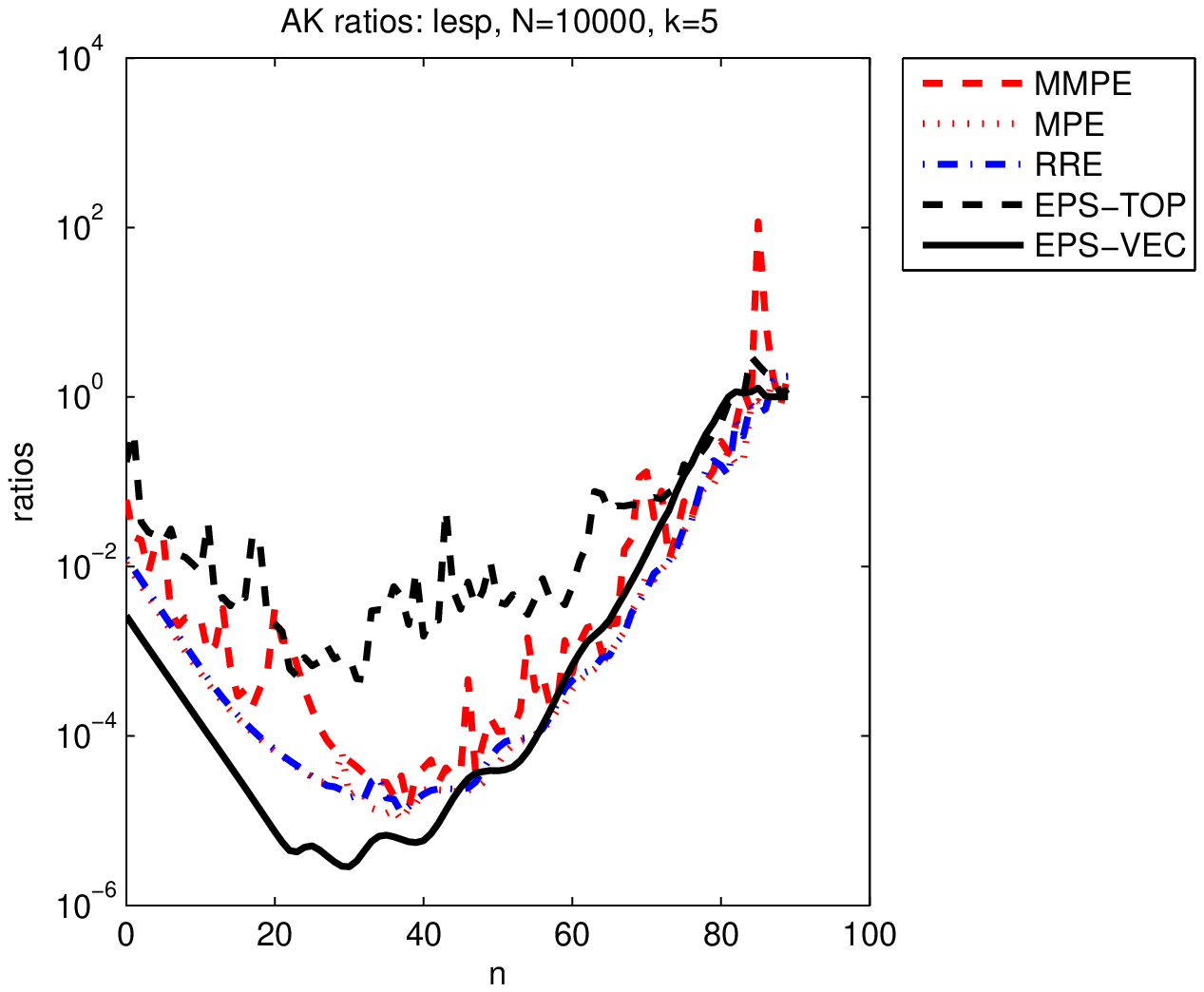}
 \end{tabular}
 \vspace{-0.7cm}
\caption{AK algorithm: errors and ratios for {\tt lesp} matrix, $N=10000$, $k=5$.}
\label{lesp}
\begin{tabular}{@{}c@{}c@{}}
 \includegraphics[width=0.5\textwidth]{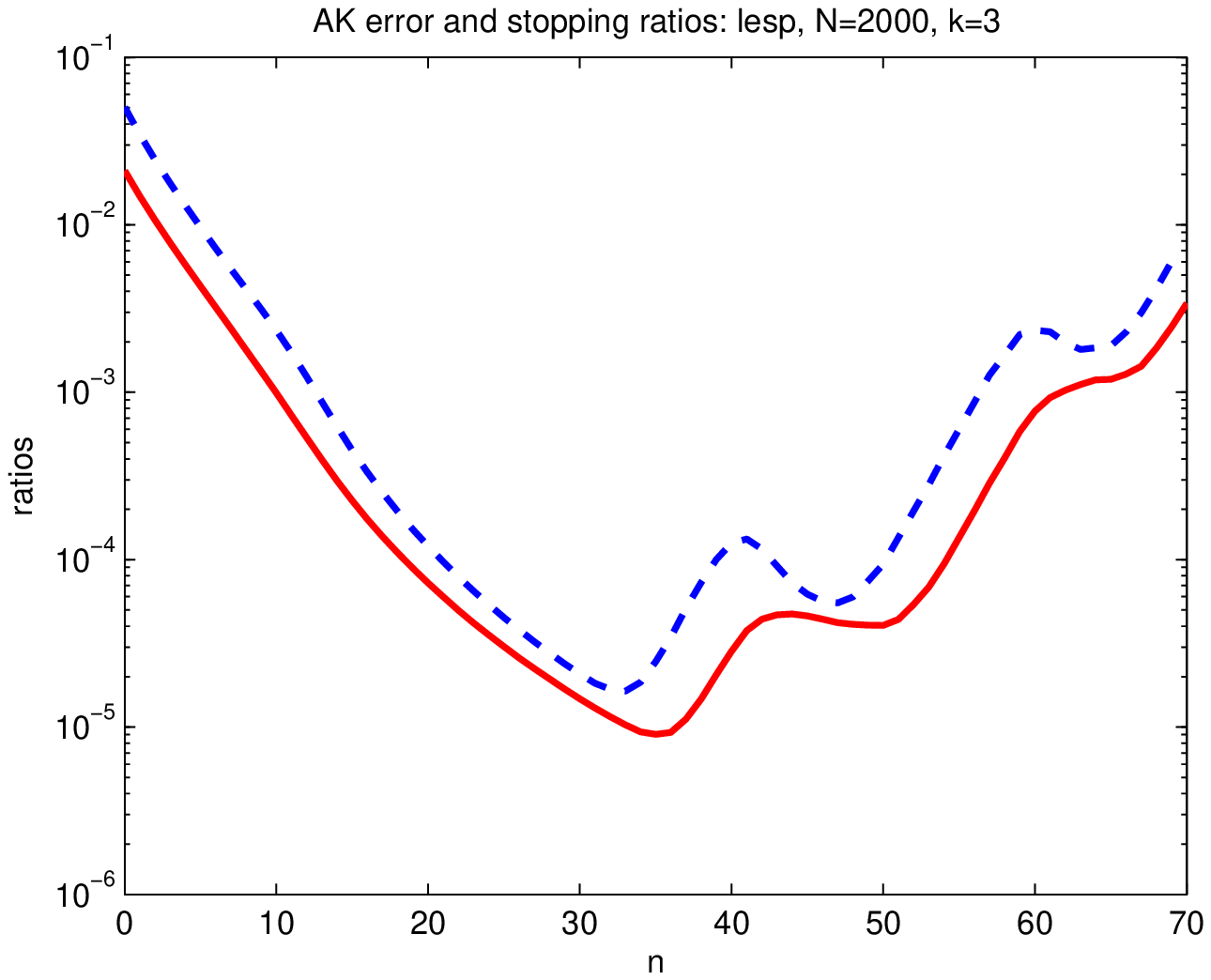}@
  \includegraphics[width=0.5\textwidth]{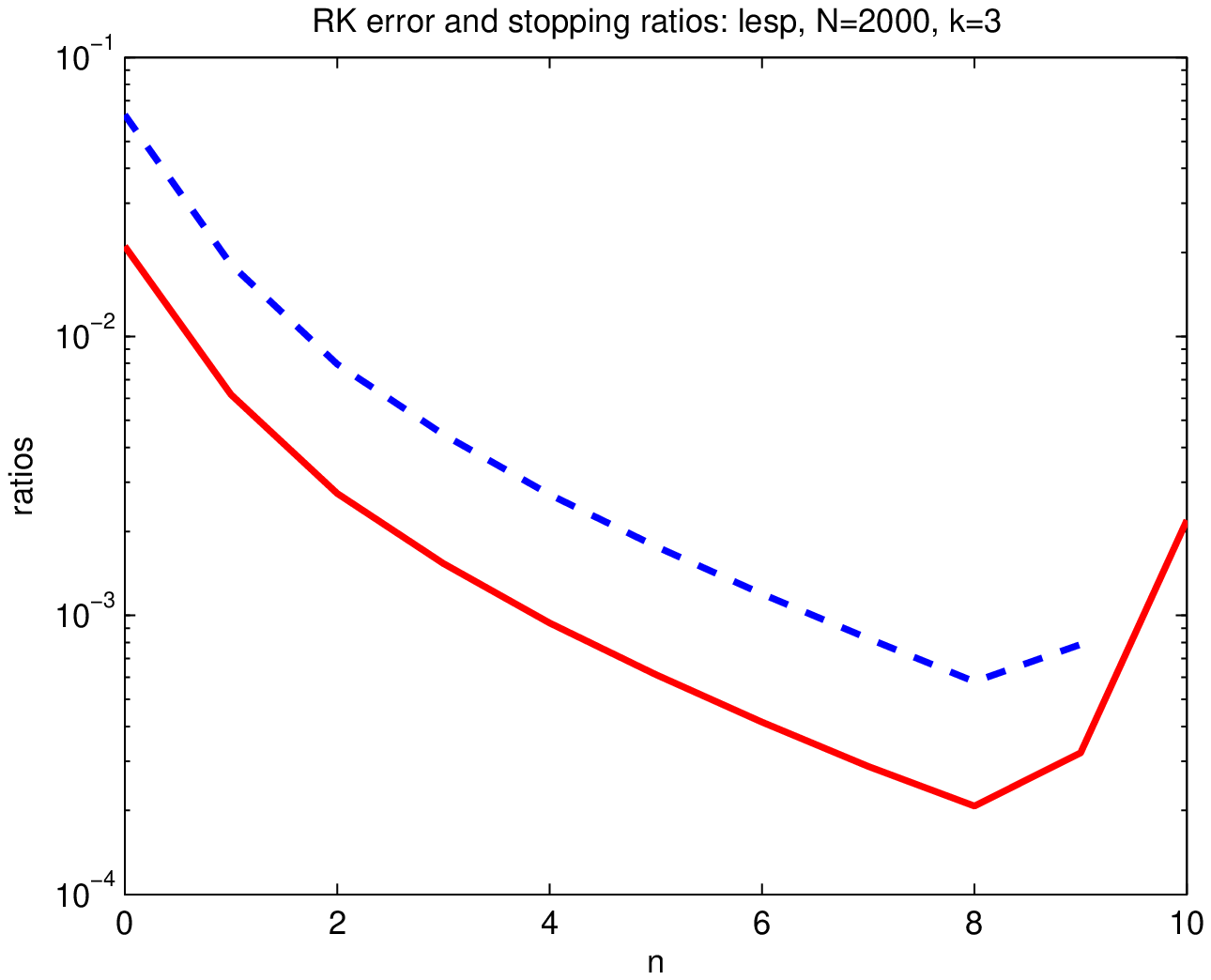}
 \end{tabular}
 \vspace{-0.7cm}
\caption{AK (left) and RK (right) algorithms: errors (solid) and stopping (dashed) ratios for {\tt lesp} matrix, $N=2000$, $k=3$.}
\label{stop}
\end{figure}

In these figures we see that all methods achieve a good precision with an advantage for the vector $\varepsilon$--algorithm. Moreover, its convergence is smoother.
The ratios grow up because all methods almost stagnate when a good precision is attained
while the error of Kaczmarz' method continues to decrease slowly.
In particular, for RK algorithm, the vector $\varepsilon$--algorithm attains its full precision after 4 iterations
while the AK algorithm needs more iterations. For this example, the dominant eigenvalue of $A$ is $0.8732178$, and the second one is $0.3170877$. Thus, according to the theoretical results of Section \ref{qqq}, a good acceleration is observed with $k=1$ for all procedures.

\vspace{0.2cm}

The Figure \ref{clementAKRK}  shows the ratios for
the {\tt clement} matrix, $\kappa(A)\simeq 1.13145 \times 10^{83}$, a tridiagonal matrix with zero diagonal entries, again with  $N=1000$ and $k=5$. For this example, the MPE and the RRE coincide. Again the vector $\varepsilon$--algorithm is the best.

\vspace{0.2cm}

For {\tt toeppen}, a pentadiagonal Toeplitz matrix, with $N=1000$, $k=5$, we have the Figure \ref{toeppen}. Notice that the MMPE and the
topological $\varepsilon$--algorithm do not work well and that this behavior could be due to the choice of the
vectors ${\mathbf{y}}$ and ${\mathbf{y}_i}$. In fact, these choices can affect the results in a quite serious way. For instance, these methods sometimes exhibit better convergence and acceleration  with ${\mathbf{y}}=(1, \ldots, 1)^T$ and ${\mathbf{y}_i}= {\mathbf{e}_i}$. With the RK algorithm, considering the first 50 iterations for $k=8$, the vector $\varepsilon$--algorithm attains a ratio
of $10^{-9}$ at iteration 20 and after 24 iterations a division by zero occurs. The ratios for the RRE and the MPE have a minimum
of $10^{-10}$ at iteration 25. The
topological $\varepsilon$--algorithm diverges from the beginning. The MMPE exhibits an erratic convergence and the ratio goes down to $10^{-7}$ at the iteration 42, and then increases.

\vspace{0.2cm}

We also tried our procedures on a bigger matrix.
The results for   {\tt lesp}, a tridiagonal matrix with real, sensitive eigenvalues, $N=10000$, $k=5$, $\kappa(A)\simeq 6.9553 \times 10^3$, with the AK algorithm are given in Figure \ref{lesp}. We see that Kaczmarz' method and the
acceleration procedures all attain full precision after 90 iterations. Thus all ratios will grow up.
 However, the vector $\varepsilon$--algorithm attains an error of less than $10^{-11}$
after about 20 iterations while  Kaczmarz's method has only an error of order $10^{-3}$ at the same iteration.

\vspace{0.2cm}

An important point in any iterative method is to have a quite reliable stopping criterion. Usually such iterations are stopped by using the residuals. However such a computation will need a matrix--vector product and, in our case, one of the interests of Kaczmarz' method will be lost. Thus,
since the results given by our acceleration procedures often stagnate when some precision is attained while those of Kaczmarz' continue to decrease, the iterations can be stopped as soon as the following ratios grow up significantly

$$\frac{\|{\Delta \mathbf{z}}_n\|}{\|{\Delta \mathbf{x}}_{n+\ell}\|}~~\mbox{(AK)} \quad \mbox{and} \quad
\frac{\|{\Delta \mathbf{z}}_n\|}{\|{\Delta \mathbf{x}}_{(n+1)(\ell+1)}\|}~~\mbox{(RK)}.$$

An example with the vector $\varepsilon$--algorithm is given in the Figure \ref{stop} (left: AK algorithm, right: RK algorithm). The solid line corresponds to the norm of the error and the dashed one to the ratio for the stopping criterion.

However, it must be noticed that, in the case of the RK algorithm, this stopping criterion involves iterates of Kaczmarz' method that have not been computed and used in the acceleration procedures. Thus, it is not usable in practice. Since
the quantities $\|{\Delta \mathbf{z}}_n\|$ usually decrease rapidly, the iterates can be stopped when it is small enough and begins to stagnate or even to grow up.

\vspace{0.2cm}

In our examples, we also add a white noise between $10^{-2}$ and $10^{-8}$ to the vector $\mathbf b$.
The norm of the error of the results obtained by our acceleration procedures
attains the level of the noise in most cases.

For the matrix {\tt baart} of the Matlab Regularization toolbox  \cite{hansenregu} of dimension 120 whose condition number is
$2.28705 \times 10^{18}$, an error of the form
$\mathbf e=\delta \|\mathbf b\| \mathbf u/ \sqrt N$ where $\delta$ is between $10^{-2}$ and $10^{-8}$ and $\mathbf u$ is a vector whose components are random variables from a normal distribution with mean 0 and standard deviation 1, was added to ${\mathbf{b}}$.
The norm of the error achieved with the vector $\varepsilon$--algorithm goes down to $10^{-0.75}$ for $k=1,2$, and $3$, which is a little bit better than the results obtained in \cite{cbmrzpn}.

Let us mention that the stopping criterion given above only works correctly for small noises.
Maybe, it because the vector $\mathbf b$ is not involved in it.

\section{Conclusions}

 From our numerical results, it seems that the vector $\varepsilon$--algorithm is the best procedure for accelerating Kaczmarz' method.
However, the recursive implementation of the other procedures has to be tested numerically to see if it leads to better and more stable results.
In our numerical examples, we tried several values of $k$. Although, when the dominant eigenvalue of $A$ is well separated, $k=1$ leads to quite a good acceleration, it seems that higher values produce better results in general. Anyway,
the choice of $k$ and of the vectors $\mathbf y$ and $\mathbf y_i$ are important points which need to be studied more deeply.
Other recursive algorithms, such as those developed by Germain--Bonne \cite{bgb}, or the VTT and the BVTT
\cite{bvt}, not considered in this work, have also to be tried.

The acceleration of the
Symmetric Kaczmarz' and the Randomized Kaczmarz' methods, which are also sequential row--action methods that update the solution using one row of $A$ at each step, of other methods
of the MAP class, of the
SIRT  (Simultaneous Iterative Reconstruction Techniques) methods, has to be considered.
Finally, applications to tomography and, in general, to image reconstruction have to be considered.



\vskip 2mm

\noindent {\bf Acknowledgments:} We would like to thank Zhong--Zhi Bai for a discussion about Meany's inequality,
Andrezj Cegielski for interesting exchanges on the terminology for the various ways of using Kaczmarz's method, and Paolo Novati for consulting about the numerical examples.
This work was partially supported by
University of Padova, Project 2010 no. CPDA104492.

\end{document}